\documentclass[10pt]{article}
\usepackage{filecontents}

\usepackage[letterpaper,text={6.7in,9in},top=0.8in,bottom=0.8in,centering]{geometry}

\usepackage{pdflscape}
\usepackage{amsmath,amssymb}
\usepackage{fleqn}
\usepackage{tikz}
\usepackage{colortbl}
\usepackage{booktabs}

\usepackage{natbib}
\usepackage[figuresright]{rotating}

\usepackage{url}
\usepackage{algorithmic}
\usepackage{algorithm}

\usepackage{amsmath,amsthm}

\theoremstyle{definition}

\theoremstyle{remark}

\title{Lot sizing problem integrated with cutting stock problem in a paper industry: a multiobjective approach}

\author{Campello B. S. C., Oliveira W. A.$^{\ast}$, Ayres A. O. C. and Ghidini C. T. L. S.}

\begin{document}
\date{}

\maketitle

\begin{center}
{\footnotesize 
*Corresponding author\\
$^1$School of Applied Sciences, University of Campinas, R. Pedro Zaccaria, 1300, Limeira 13484-350, S\~ao Paulo, Brazil\\
E-mails: betacampello@gmail.com [Campello] / washington.oliveira@fca.unicamp.br [Oliveira] / amandaortega@gmail.com [Ayres] / carla.ghidini@fca.unicamp.br [Ghidini]
}\end{center}

\begin{abstract}
In this work, we use a multiobjective approach to address the lot sizing problem integrated with the cutting stock problem in a paper industry. We analyze the trade-offs and correlations which exist among the costs and their decision variables. Considering some of our computational results, if we decrease the production costs, then we increase the waste of material of the cutting process and vice versa. Thereby we show the importance of the multiobjective approach in allowing multiple answers to the decision maker, using Pareto optimal solutions set. Several tests were performed to check the quality of our approach.
\end{abstract}


\section{Introduction}\label{sec:introduction}
In a paper industry, the primary manufacturing process consists of producing material objects according to demanded lengths and paper types, setup and capacity of processing machines and the extra of these objects can be stored to use in the long-term planning horizon. Lot sizing problem (LSP) consists of planning the quantities of material objects that will be produced in a period to minimize the production, storage and setup costs, without exceeding the capacity of the machine and meeting the consumer demands. We can find two surveys,~\cite{drexl1997lot} and~\cite{jans2008modeling}, regarding relevant papers which cover this subject.

Following the objects are cut into smaller pieces. In the one-dimensional case, {\it cutting pattern} is the configuration in which the piece lengths are organized into the object length to be cut, whereas it can generate a waste of material. These cut pieces can be stored to meet future demands in the long-term planning horizon. Cutting stock problem (CSP) consists of choosing the best cutting patterns and deciding their quantities to be used to minimize waste of material and storage costs and meet the customer demands. CSP has been widely studied in the literature, part of these studies can be found in three surveys (\cite{sweeney1992cutting,cheng1994cutting,de2002lp}) dedicated to relevant papers about CSP.

There is a dependency between these two problems. Since in CSP the waste of material depends on the length and the amount of produced objects in LSP. The latter should consider the cutting process in the production planning of objects. Otherwise, it can generate far more waste of material in CSP (\cite{poltroniere2008,gramani2009lagrangian}). From this fact, some recent researches proposed to integrate LSP with CSP in a mono-objective approach, with the purpose of improving the overall costs of production.

In a study about integrated process optimization,~\cite{arbib2005integrating} integrate hierarchical decision levels and functional areas (short-term operations, mid-term planning, production and purchase of materials), and by using leftovers of material they achieve on the average a further reduction of about 43\% concerning the overall costs. With some similarities to the model discussed here, in a furniture industry case (\cite{gramani2006combined}), the overall costs decrease 13\% in comparison to the disconnected approach. Another example occurs in a paper industry (\cite{malik2009integrated}), where it reduces the overall costs in about 25\%. Similar mono-objective approaches are discussed by~\cite{hendry1996cutting, nonaas2000combined, ghidini2007solving, poltroniere2008, suliman2012algorithm} and~\cite{silva2014integrating}.

The mono-objective integrated approach provides just only one optimal solution to be examined by a decision maker, which makes it difficult to assess their effectiveness, because of the trade-off between LSP and CSP. Instead, the multiobjective integrated approach allows a plethora of solutions (\cite{miettinen2012nonlinear}), allowing a more extended study about the trade-off between these two problems, by considering the features among different optimal solutions. Thus, we can infer that the multiobjective integrated approach is more convenient for the manufacturing industries because they would organize the production under their own criteria from different optimal solutions.

We did not find in the literature multiobjective integrated approaches for these two problems, but only texts where they are considered isolated. In CSP case,~\cite{wascher1990lp} minimizes, from a multiobjective approach, the storage costs, the extra cost of inputs and the waste of material during the cutting process. And~\cite{golfeto2009genetic} propose a bi-objective model to minimize the waste of material and the machine setup, and they use a symbiotic genetic algorithm to solve their proposed model. In LSP case,~\cite{ustun2008integrated} connect the selection of raw material suppliers with LSP costs. And~\cite{rezaei2011multi} introduce two multi-objective mixed integer non-linear models for LSP multi-period involving multiple products and multiple suppliers, where the cost, the quality level, and the service level are considered in the objective vector, and after applying a genetic algorithm, the total costs are reduced significantly.

The central aspect of this contribution is to evaluate the trade-off between LSP and CSP by analyzing the cost variations of these two problems when we simultaneously minimize them. Also, we separated the objective function in some other functions to compare, detect, analyze and study potential trade-offs and correlations at the heart of the problem. In this work, we show that there are several solution options for the integrated problem. We obtained these solutions by applying two multiobjective traditional methods of resolution. Last, we analyze the correlations among some decision variables and present the advantages of studying LSP integrated with CSP using a multiobjective approach. 

The organization of the paper is as follows. Section~\ref{sec:modeling} gives a formal description of the integrated problem as a multiobjective optimization programming. Section~\ref{sec:methods} describes briefly the two resolution methods used. In Section~\ref{sec:numericalexperiments}, we present the main parameters used in the computational experiments and examine the computational results. In the end, Section~\ref{sec:conclusions} presents the conclusions.

\section{Mathematical modeling}\label{sec:modeling}

\subsection{Description of the problem}

\begin{figure}[h]
\begin{center}
\includegraphics[width=7cm]{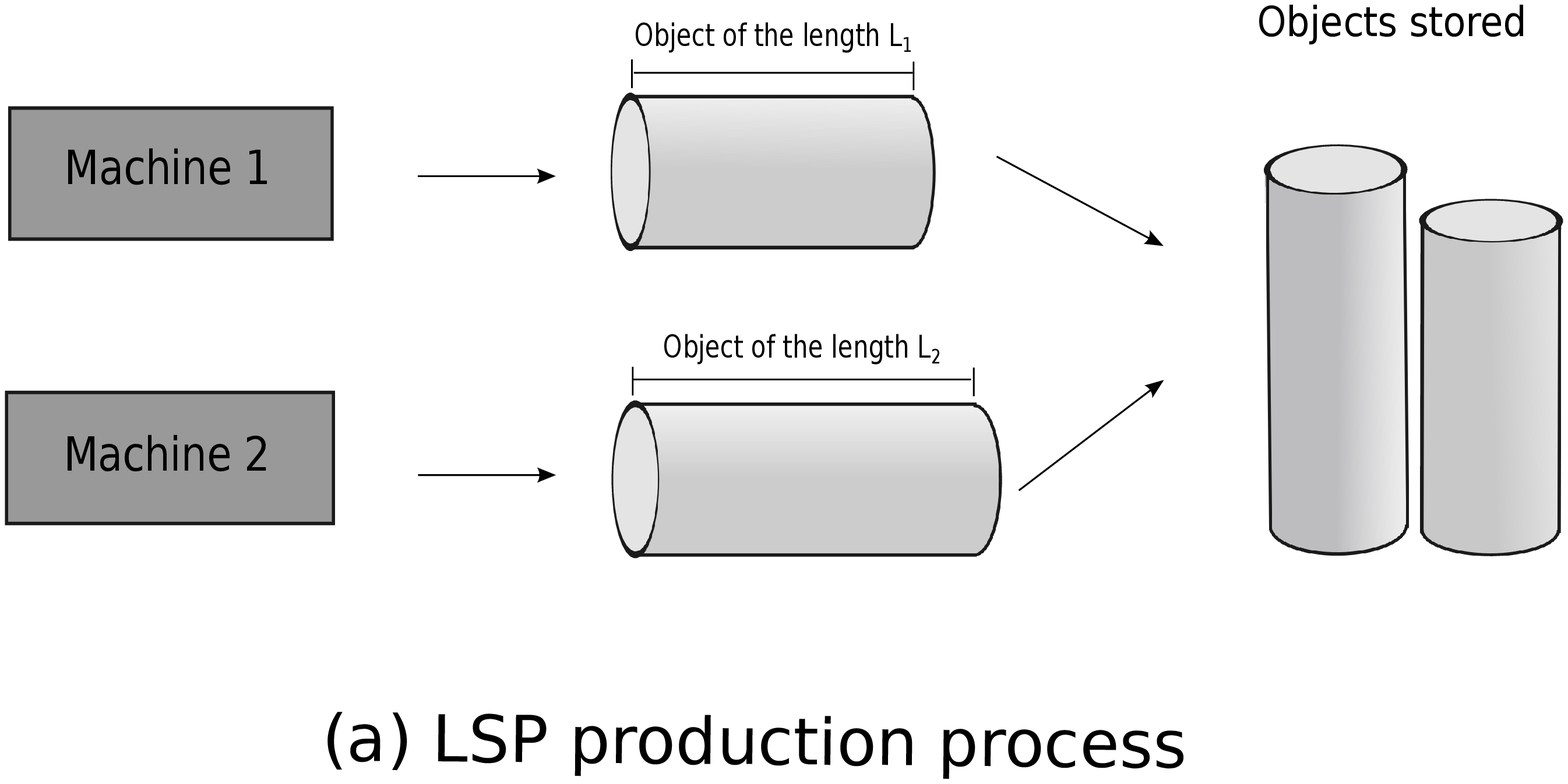}\hspace{1.4cm}
\includegraphics[width=7cm]{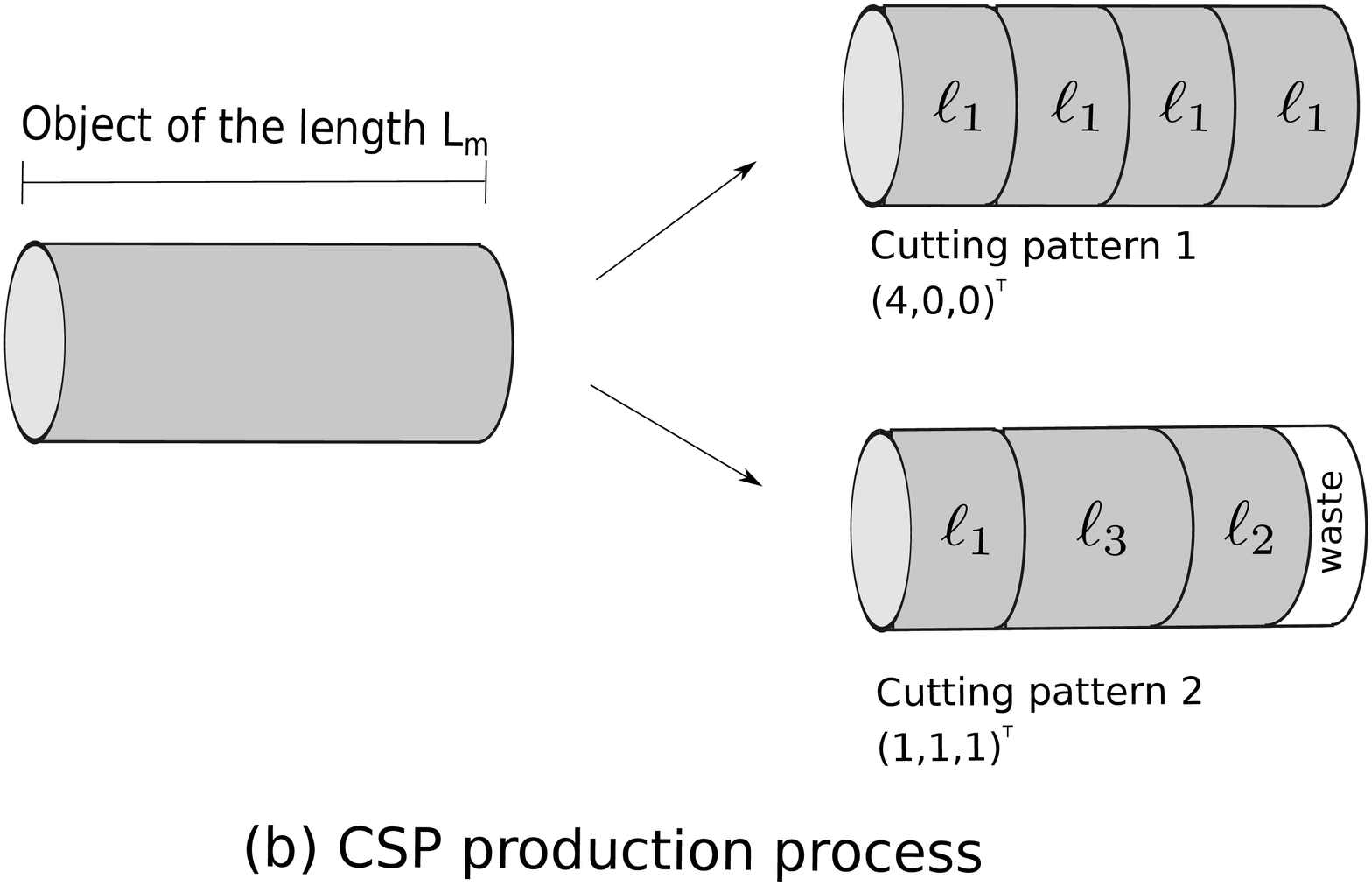}
\caption{Papermaking processes}
\label{fig:producaoobjeto}
\end{center}
\end{figure}

The first part of the papermaking processes consists of transforming the cellulose pulp into a continuous and smooth sheet, which is rolled up as finished paper. The results of this process are some material objects (jumbos, master rolls, bars, or simply objects) that can be of different lengths. We assume $m$ available machines, $m=1,\ldots,M$, producing objects of length $L_m$. There are inventory costs if the object is stored, and setup costs if the machine $m$ is used, because of the waste of material during the grammage exchange. Besides, each machine has a limited production capacity. Figure~\ref{fig:producaoobjeto}(a) illustrates the LSP manufacturing process, where two machines produce two objects of lengths $L_1$ and $L_2$, which are stored. Therefore, in this work, LSP consists of deciding whether a machine must be used, the amount of objects to be produced in a machine, and the quantities of objects to be stored to minimize production, storage and setup costs. 

The second stage of the papermaking process consists of cutting each object into smaller parts (smaller jumbos, smaller rolls, items, formats, or only pieces). Each demanded piece $i$, $i=1,\ldots,Nf$, is of length $\ell_i$. The pieces can be stored to meet future demands in the planning horizon, but this causes storage costs. Figure~\ref{fig:producaoobjeto}(b) illustrates this second papermaking process, where a waste of material is generated according to the chosen cutting pattern. Hence, in this work, CSP consists of deciding the frequencies to cut the object according to each cutting pattern and the amount of pieces needed to be stored for the purpose of minimizing the waste of material costs and the storage costs.

Here, we connect these two papermaking processes as a unique problem and study its solutions from a multiobjective optimization approach, where the decision variables are linked, and the manufacturing processes impact each other.

\subsection{Multiobjective model}

We introduce our mathematical formulation for modeling the production process of a paper industry as a bi-objective optimization problem. We used a set of constraints similar to the one presented in~\cite{poltroniere2008}. From the same study, we turned its mono-objective function into two new objective functions, one for LSP and another one for CSP.  We consider the following notation.\\   

\noindent{\bf\textit{Index}:}

\noindent \textbf{{LSP:}}

\noindent $T$: period quantities in the planning horizon; period index is $t = 1,\ldots, T$;

\noindent $K$: paper grammage type quantities; grammage index is $k = 1,\ldots, K$;

\noindent $M$: machine type quantities; machine index is $m = 1,\ldots, M$; each machine type $m$ produces objects of length $L_m$.

\noindent \textbf{{CSP:}}

\noindent $N_m$: cutting pattern type quantities of the object of length $L_m$; cutting pattern index is $j=1,\ldots,N_m$; 

\noindent $Nf$: demanded piece type quantities; piece index is $i=1,\ldots,Nf$; 

\noindent $\{1,\ldots,Nf\}\equiv S(1)\cup S(2)\cup...\cup S(K)$, where $S(k)$=\{$i|$ the piece type $i$ is of grammage $k$\}.\\

\noindent{\bf\textit{Parameters}:}

\noindent \textbf{{LSP:}}

\noindent $c_{kmt}$: cost to produce the object of grammage $k$ in the machine type $m$ in the period $t$;

\noindent $h_{kt}$: cost to store the object of grammage $k$ at the end of the period $t$;

\noindent $s_{kmt}$: setup cost to produce the object of grammage $k$ in the machine type $m$ in the period $t$;

\noindent $C_{mt}$: capacity of the machine type $m$ in the period $t$;

\noindent $\rho_k$: specific weight of objects of grammage $k$;

\noindent $D_{kt}$: demand of paper of grammage $k$ in the period $t$;

\noindent $b_{km}$: object weight of grammage $k$ produced in the machine type $m$ ($b_{km}=\rho_k L_m$); 

\noindent $f_{km}$: amount of waste of material of grammage $k$ during the production of an object in the machine type $m$. 

\noindent \textbf{{CSP:}}

\noindent $\sigma_{it}$: cost to store the piece type $i$ in the period $t$; 

\noindent $cp_{kt}$: cost of waste of material of grammage $k$ in the cutting process in the period $t$;

\noindent ${\bf d}_{kt}$: vector of demanded pieces; each component $d_{ikt}$ represents the amount of demanded piece type $i$, $i \in S(k)$, in the period $t$; 

\noindent $\eta_{ik}$: weight of the piece $i$ of grammage $k$ ($\eta_{ik}=\rho_k \ell_i$); 

\noindent ${\bf a}_{jm}$: vector of cutting patterns; each component $a_{ijm}$ represents the amount of demanded piece type $i$ in the cutting pattern $j$ of the object of length $L_m$; 

\noindent $p_{jm}$: amount of waste of material in the cutting pattern $j$ of the object of length $L_m$;

\noindent $Q$: a sufficiently large number.\\

\noindent {\bf\textit{Decision variables}:}

\noindent \textbf{{LSP:}}

\noindent $x_{kmt}$: number of objects of grammage $k$ produced in the machine type $m$ in the period $t$; 

\noindent $w_{kmt}$: number of objects of grammage $k$ produced in the machine type $m$ and stored at the end of the period $t$; 

\noindent $z_{kmt}$: binary variable; it takes the value 1 if the object of grammage $k$ is produced in the machine type $m$ in the period $t$, and 0 otherwise.\\\\

\noindent \textbf{{CSP:}}

\noindent $y^j_{kmt}$: number of objects of grammage $k$ produced in the machine type $m$  in the period $t$, which are cut according to the cutting pattern $j$; 

\noindent ${\bf e}_{kt}$: vector of storage pieces; each component $e_{ikt}$ represents the amount of piece type $i$ of grammage $k$ stored at the end of the period $t$.\\

The LSP total costs are modeled by  function $F_1(x,w,z)$. It measures simultaneously three primary goals to be minimized, the production, the storage of object and the machine setup costs.
\begin{eqnarray}\label{eq:f1}
\mathrm\ F_1(x,w,z) \ = \ \sum_{t=1}^{T} \sum_{m=1}^{M} \sum_{k=1}^{K} (c_{kmt} x_{kmt} + h_{kt} b_{km} w_{kmt} + s_{kmt} z_{kmt})
\end{eqnarray} 

The CSP total costs are modeled by function $F_2(y,e)$. It measures simultaneously two primary goals to be minimized, the waste of material and the storage of piece costs.
\begin{eqnarray}\label{eq:f2}
\mathrm\ F_2(y,e) \ = \ \sum_{t=1}^{T} \sum_{k=1}^{K} cp_{kt} \sum_{m=1}^{M} \sum_{j=1}^{N_m} p_{jm} y^{j}_{kmt} + \sum_{t=1}^{T} \sum_{k=1}^{K} \sum_{i \in S(k)} \sigma_{it} \eta_{ik} e_{ikt}  
\end{eqnarray} 

\noindent{\bf\textit{Bi-objective mathematical model}:}

\begin{align}
\text{Min }   & F \ = (F_1(x,w,z),F_2(y,e))\label{eq:f}\\ 
\text{s.t.\ \textcolor{white}{.}  }&  \nonumber\\
&\sum_{m=1}^{M} (b_{km} x_{kmt} + b_{km} w_{k,m,t - 1} - b_{km} w_{kmt}) \geq D_{kt}, \quad k = 1,\ldots, K,\ t = 1,\ldots,T;\label{eq:producao_objeto}\\
&\sum_{k=1}^{K} (b_{km} x_{kmt} + f_{km} z_{kmt}) \leq C_{mt}, \quad m = 1,\ldots, M,\ t = 1,\ldots, T;\label{eq:capacidade_maquina}\\
&x_{kmt} \leq Q z_{kmt}, \quad k = 1, ..., K,\ m = 1,\ldots, M,\ t = 1,\ldots, T;\label{eq:setup}\\
&\sum_{m = 1}^{M} \sum_{j = 1}^{N_m} {\bf a}_{jm} y^{j}_{kmt} + {\bf e}_{k, t - 1} - {\bf e}_{kt} = {\bf d}_{kt}, \quad k = 1,\ldots, K,\ t = 1,\ldots, T;\label{eq:corte_itens}\\
&\sum_{j = 1}^{N_m} y^{j}_{kmt} = x_{kmt} + w_{k,m,t - 1} - w_{kmt}, \quad k = 1,\ldots, K,\ m = 1,\ldots, M,\ t = 1,\ldots, T;\label{eq:acoplamento}\\
&w_{km0} = 0,\ e_{k0} = 0, \quad k = 1,\ldots, K,\ m = 1,\ldots, M;\label{eq:condicao_inicial}\\
&x_{kmt},\ w_{kmt} \in \mathbb{Z}^+,\ z_{kmt}\in \{0, 1\},\quad k = 1,\ldots, K,\ m = 1,\ldots,M,\ t = 1,\ldots, T;\label{eq:condicao_xwz}\\
&y^{j}_{kmt},\ e_{kt} \in \mathbb{Z}^+,\quad j = 1,\ldots ,N_{m},\ k = 1,\ldots, K,\ m = 1,\ldots ,M,\ t = 1,\ldots,T.\label{eq:condicao_ye}
\end{align}

By using the index, parameters, and variables described above, we explained the constraints set as follows. Constraints~(\ref{eq:producao_objeto})--(\ref{eq:setup}) are related to LSP. The meeting of paper demand, measured in total weight, is guaranteed by~(\ref{eq:producao_objeto}). As described by~\cite{poltroniere2008}, we can consider the total paper demand $D_{kt}$ as being the demand of pieces plus the waste of material from cutting process (both in tonnes), it means, $D_{kt}=\sum_{i \in S(k)} \eta_{ik} d_{ikt} + waste$. But, in this work, we define $waste$ equal to zero, to facilitate the modeling of this demand constraint and let $D_{kt}$ to be a parameter, and we allow the total weight of produced objects is equal to or greater than the total paper demand. Constraints~(\ref{eq:capacidade_maquina}) make sure that machine capacity is not exceeded. If a machine is used,~(\ref{eq:setup}) express the fact that its setup costs must be computed. Constraints~(\ref{eq:corte_itens}) are related to CSP, they imply that the demand of pieces is met in each period, and for each $k$ and $t$, we note that ${\bf a}_{jm}$, ${\bf e}_{kt}$ and ${\bf d}_{kt}$ have the same dimension $|S(k)|$. The integration between LSP and CSP is enforced by~(\ref{eq:acoplamento}), the amount of produced objects plus the stock balance are equal to the cut object quantities. Constraints~(\ref{eq:condicao_inicial}) define the opening stock levels and the last sets of constraints,~(\ref{eq:condicao_xwz}) and~(\ref{eq:condicao_ye}), simply specify the allowed values for the variables.

\section{Methods of resolution}\label{sec:methods}

Two well-known methods were used to obtain Pareto optimal solutions. Taking into account~\cite{miettinen2012nonlinear}, we briefly describe these methods.\\

\noindent{\textit{Definition}}: A feasible point $u^*=(x^*,y^*,z^*,w^*,e^*)$ is said to be a Pareto optimal solution of the problem~(\ref{eq:f})--(\ref{eq:condicao_ye}), if there is not another feasible point $u=(x,y,z,w,e)$ such that $(F_1(u),F_2(u))\leq (F_1(u^*), F_2(u^*))$ and $(F_1(u),F_2(u))\neq (F_1(u^*), F_2(u^*))$.\\

\noindent{\bf\textit{Weighting method}}

In the weighting method, we associate each objective function with a weighting coefficient and minimize the weighted sum of the objective functions (i.e., a weighted mono-objective problem). For each iteration of the computational tests, we change the weighting coefficients to obtain several Pareto optimal solutions. In the objective function ($p_1 F_1 + p_2 F_2$), we use weights $p_i$, in the interval (0,1), satisfying $p_1+p_2=1$.

In practice, the objective functions need to be normalized to avoid that one of them prevails upon the other, because of its greater magnitude. We use the normalization,
\begin{align}
\mathrm \ f_{k}^{norm}(x)\ = \frac{f_{k}(x) - z_{k}^{*}}{z_{k}^{nad} - z_{k}^{*} }, k = 1,\ldots, K,
\end{align}

\noindent where $z_{k}^{*}$ is the $k$-th component of the Ideal-solution vector and $z_{k}^{nad}$ is the $k$-th component of the Nadir-solution vector (\cite{miettinen2012nonlinear}). Therefore, $f_{k}^{norm}$ is the $k$-th normalized objective function in the interval (0,1).\\

\noindent{\bf\textit{$\varepsilon$-Constraint method}}

In the $\varepsilon$-constraint method, we select one of the objective functions to be minimized, and we convert all the other ones into constraints by setting an upper bound $\varepsilon$ to each of them. If we set up the upper bounds in an appropriate manner, it is possible to reach all the Pareto optimal solutions. In this work the LSP function $F_1(x,w,z)$ is minimized and the CSP function $F_2(y,e)$ is converted into constraint. Thus,~(\ref{eq:f}) becomes 
\begin{align}
\text{Min  }   &F \ = F_1(x,w,z),\label{eq:pdlepsilon} 
\end{align}  

\noindent and we include~(\ref{eq:pceepsilon}) in the constraint set.
\begin{align}
 \sum_{t=1}^{T} \sum_{k=1}^{K} cp_{kt} \sum_{m=1}^{M} \sum_{j=1}^{N_m} p_{jm} y^{j}_{kmt} + \sum_{t=1}^{T} \sum_{k=1}^{K} \sum_{i \in S(k)} \sigma_{it} \eta_{ik} e_{ikt} \leq \varepsilon_2. \label{eq:pceepsilon}
\end{align}   

Considering the Ideal and Nadir vectors, we use upper bounds $\varepsilon_2$ by considering an interval range from Ideal to Nadir solution of CSP component. As $\varepsilon_2$ increasing from Ideal solution, it makes possible to find better solutions for LSP objective function. 

\section{Numerical experiments}\label{sec:numericalexperiments}
In this section, we measure the quality and the performance of our integrated multiobjective approach. We have solved the proposed mathematical model by using the solver CPLEX 12.6.1. The two methods of resolution were coded in Python language and executed on an Intel Quad Core i5, 24 GB of RAM and Windows operating system. We used as linking code between CPLEX and Python the \textit{docplex} library.

We use the benchmark suite of 12 instance classes of the proposed model as described in Table~\ref{tab:classes} to test our approach with the same feature found in \cite{poltroniere2008}. The full characterization of an instance involves the following parameters: the type of piece $i$ and its quantity $N_f$; the period $t$ and its quantity $T$; the type of machine $m$ and its quantity $M=2$. The lengths of the objects are $L_1=540$ and $L_2=460$ centimeters; the specific weight $\rho_k$ of all objects is 2kg/cm; and the amount of grammage is $k=1$. We present the range for each parameter described in Section~\ref{sec:modeling}, the values for the parameters were generated randomly in the interval range described below, where $h_{kt}$ and $\sigma_{it}$ are costs measured per tonne, $cp_{kt}$ is a cost measured per centimeter, $C_{mt}$, $D_{kt}$ and $f_{km}$ are measured in tonnes, $p_{jm}$ is measured in centimeter, and $c_{kmt}$ and $s_{kmt}$ are proportional to the total weight of the object.\\

\begin{table}[t]
\caption{Data of instances}
\begin{center}
\begin{tabular}{crrrrrrrrrrrr}
\toprule
Class & 1 & 2 & 3 & 4 & 5 & 6 & 7 & 8 & 9 & 10 & 11 & 12 \\ \midrule
$N_f$& 3 & 3 & 4 & 4 & 5 & 5 & 6 & 6 & 7 & 7 & 8 & 8 \\
$T$ & 3 & 4 & 3 & 4 & 3 & 4 & 3 & 4 & 3 & 4 & 3 & 4 \\ \bottomrule
\end{tabular}
\end{center}
\label{tab:classes}
\end{table}

\begin{description}
\item $c_{kmt}\in [0.015, 0.025]*b_{km}$, where $b_{km} = \rho_k L_m$;
\item $s_{kmt}\in [0.03, 0.05]*c_{kmt}$;\ $h_{kt} \in [0.0000075, 0.0000125]$;\ $f_{km} \in [0.01, 0.05]*b_{km}$; 
\item $cp_{kt} = \displaystyle \frac{\sum_{m=1}^{M}c_{kmt}}{M}*10$;\ $\sigma_{it} = 0.5*h_{kt}$;
\item $\ell_i\in [0.1, 0.3]*\displaystyle\frac{\sum_{m=1}^{M} L_m}{M}$;\ $d_{ikt}\in [0, 300]$;\ $\eta_{ik}=\rho_k \ell_i$;
\item $C_{mt} = \displaystyle\frac{b_{km}}{\sum_{m = 1}^{M}b_{km}}$*$Cap$; where $Cap$ = $\phi* \displaystyle\frac{\sum_{t=1}^{T} \sum_{m=1}^{M} \sum_{k=1}^{K} ((D_{kt}/M) + f_{km})}{M*T}$, and $\phi = 1.24$.
\end{description}

For each class, we solved twenty instances. For each instance and each method of resolution, we have chosen to use fifty weighting coefficients $p_i\in (0,1)$ in the weighting method, and fifty upper bounds $\varepsilon_2$ by considering an interval range from Ideal to Nadir solution of CSP in the $\varepsilon$-constraint method. We treated the computational experiments through three groups of tests. For each group, we tested the two methods of resolution described in Section~\ref{sec:methods} as follows.

{\bf Group 1}. \textit{Suboptimal solutions from heuristic}: It is known that the two problems in the integrated mathematical model are in general NP-hard. Thus we apply a heuristic approach to solve it in a reasonable computational time. Next, we describe the steps of the heuristic approach.

{\bf Step 1} - For each $m$, generate all the possible cutting patterns.

{\bf Step 2} - For each $m$, sort the cutting patterns in ascending order according to its waste of material and select the first $n$ cutting patterns.

{\bf Step 3} - Solve the mathematical modeling described in Section~\ref{sec:modeling} by using just the $n M$ cutting patterns selected in Step 2.

In this group of tests, we define $n=15$, in every period we use the same $n M$ cutting patterns, and we generate instances from Class 1 to Class 8 of Table~\ref{tab:classes}.\\

{\bf Group 2}. \textit{Optimal Integer solutions from Solver}: It is possible to solve some small dimensional instances directly by using CPLEX solver. Consequently, we can validate the trade-off between the objective functions without the influence of suboptimal solutions from heuristic. In this group of experiments, we generate instances from Class 1 to Class 4 of Table~\ref{tab:classes}.\\

{\bf Group 3}. \textit{Relaxing integrality of the variables}: In particular, due to the difficulty in resolving significant dimensional instances, the computational tests on the previous groups cover small and medium size instances. With the purpose of validating our proposed approach to solve large size instances, we relax the integrality of the variables, except for the binary variable $z$ regarding the setup costs. In this group of tests, we generate instances for all classes of Table~\ref{tab:classes}.\\

In the next subsection, we close the numerical experiments with statistics summarizing the test results and comments on the performances of our approach.

\subsection{Global analysis of the results}

\begin{figure}[t]
\begin{center}
\includegraphics[width=8cm]{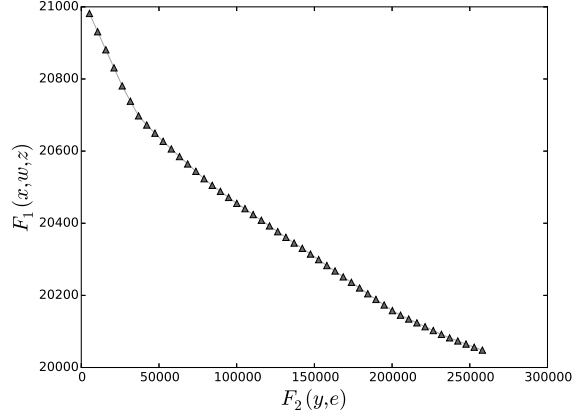}
\caption{Pareto front from the $\varepsilon$-constraint method} \label{fig:curvaParetoGrupo3Erestrito}
\end{center}
\end{figure}

The main purpose of this study is to evaluate the trade-off between LSP and CSP by modeling them as a multiobjective integrated problem. Also, from a suitable manner, we separated the objective function in some other functions to be compared. Thus we can detect, analyze and study the potential trade-offs/correlations at the heart of the problem by using the variables that have the greatest effect each other. We examine the following issues: if the increase in the CSP total costs involves the decrease in the LSP total costs or vice versa; if there is just one optimal solution for the integrated problem, i.e., if there is no trade-off; and what the major advantages achieved by using the multiobjective integrated approach are.

\begin{figure}[t]
\begin{center}
\includegraphics[width=5cm]{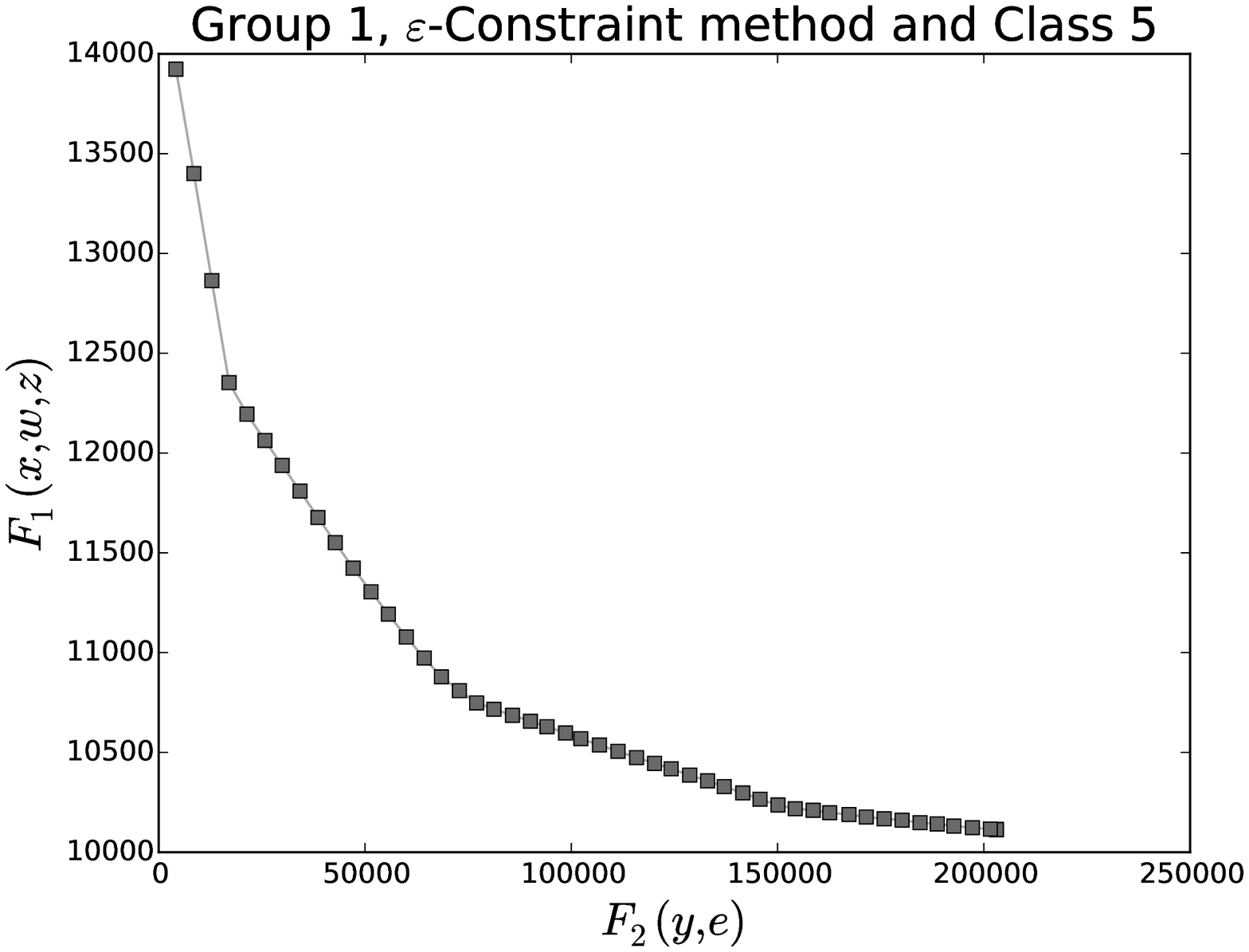} 
\includegraphics[width=5cm]{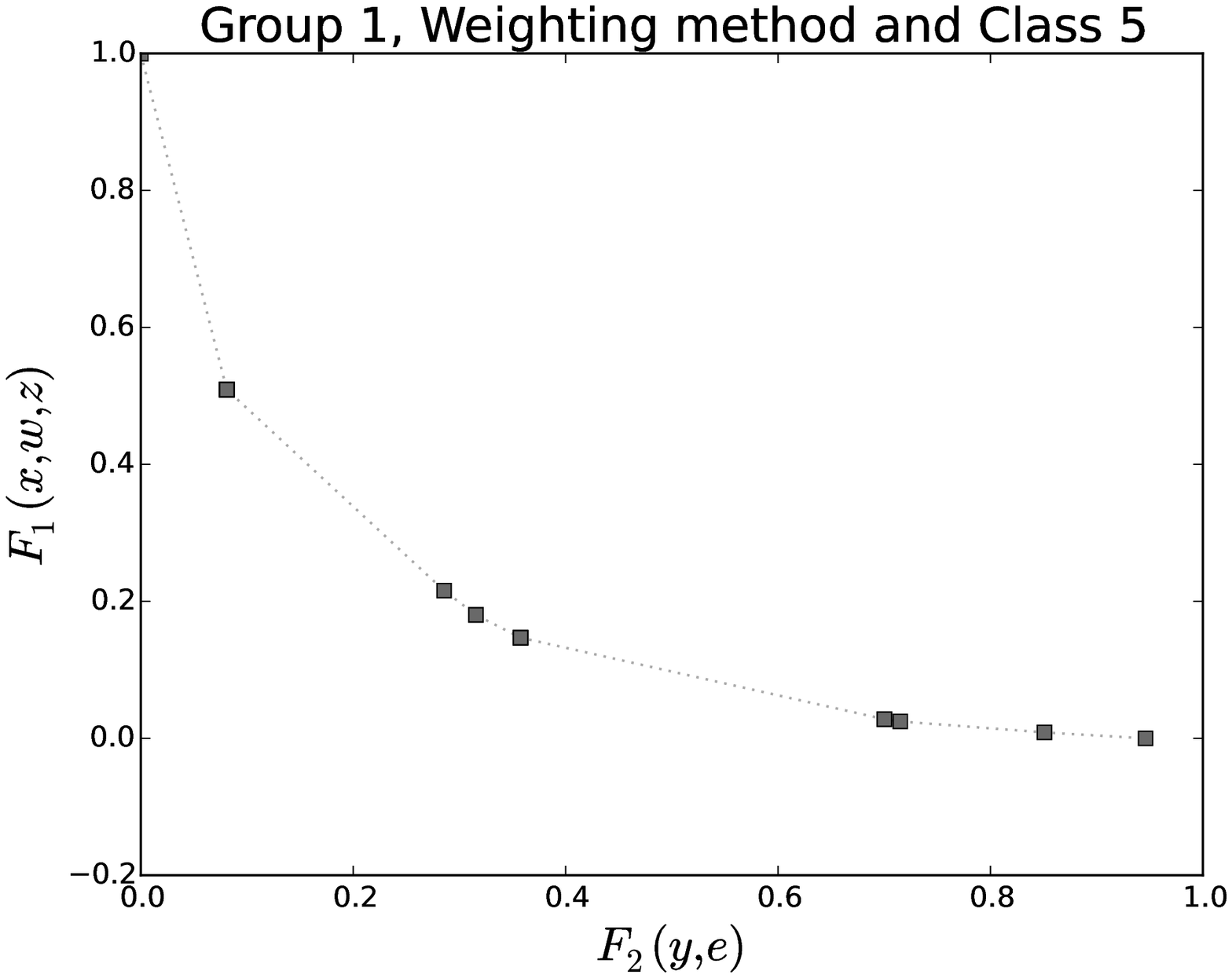} 
\includegraphics[width=5cm]{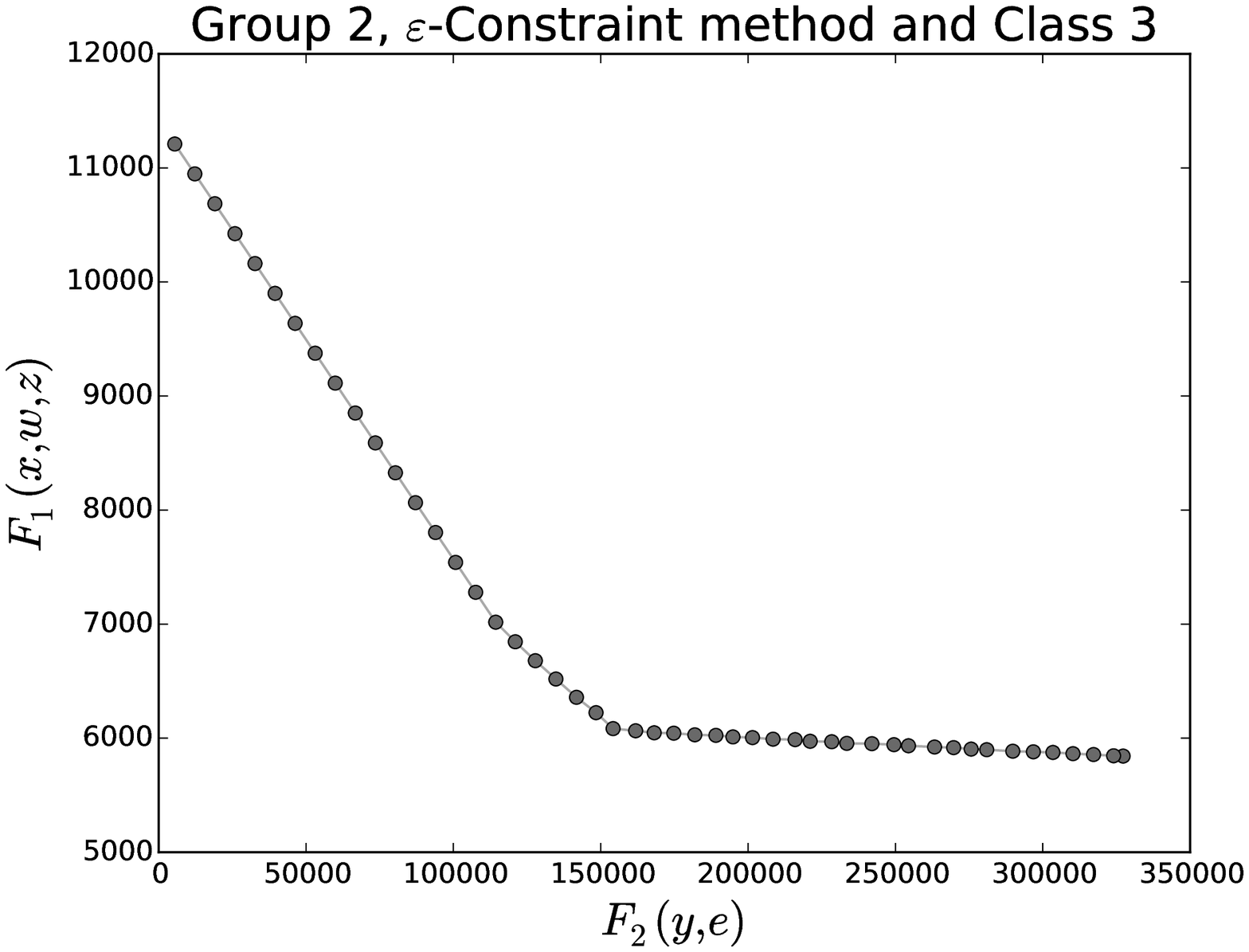}
\includegraphics[width=5cm]{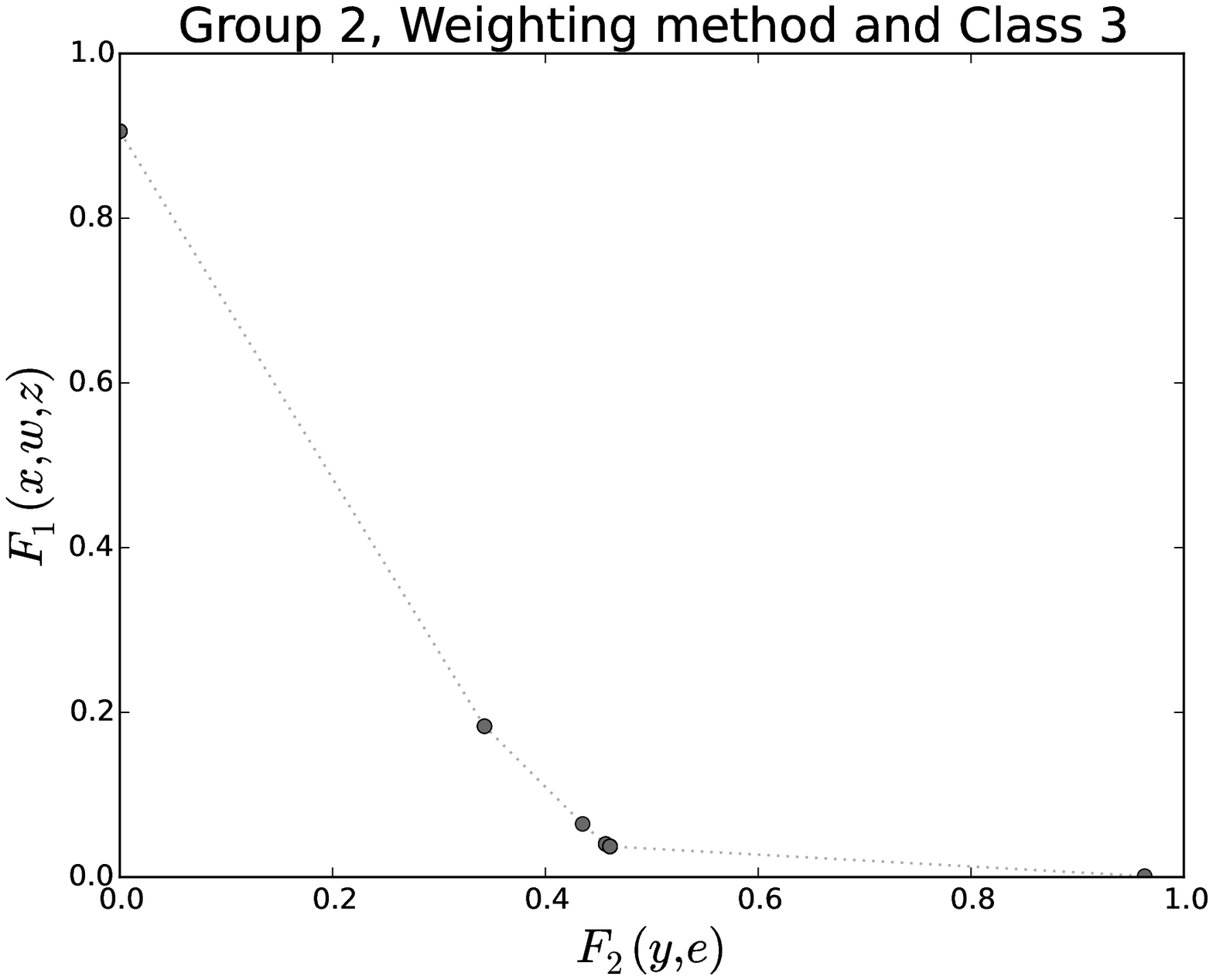}
\includegraphics[width=5cm]{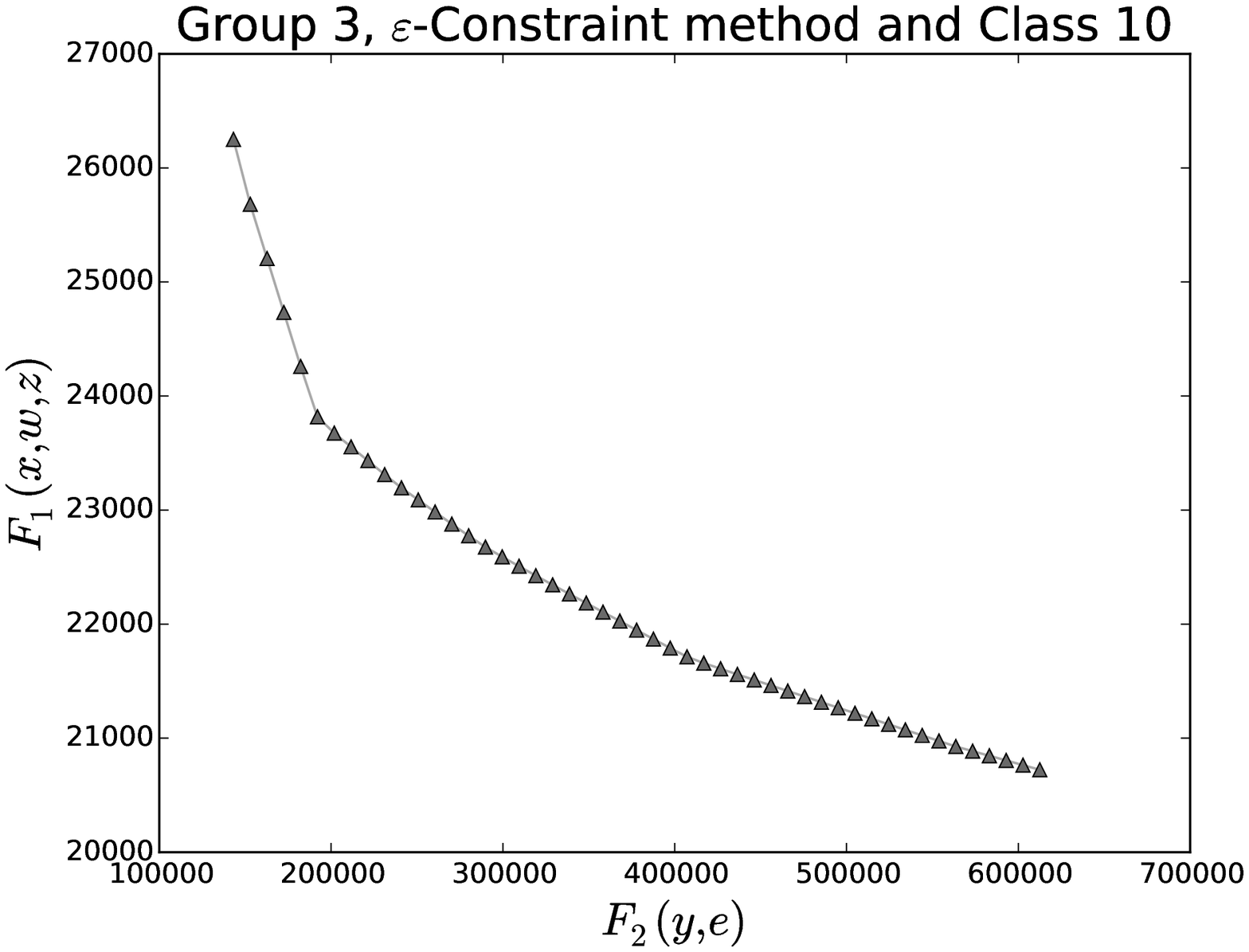}
\includegraphics[width=5cm]{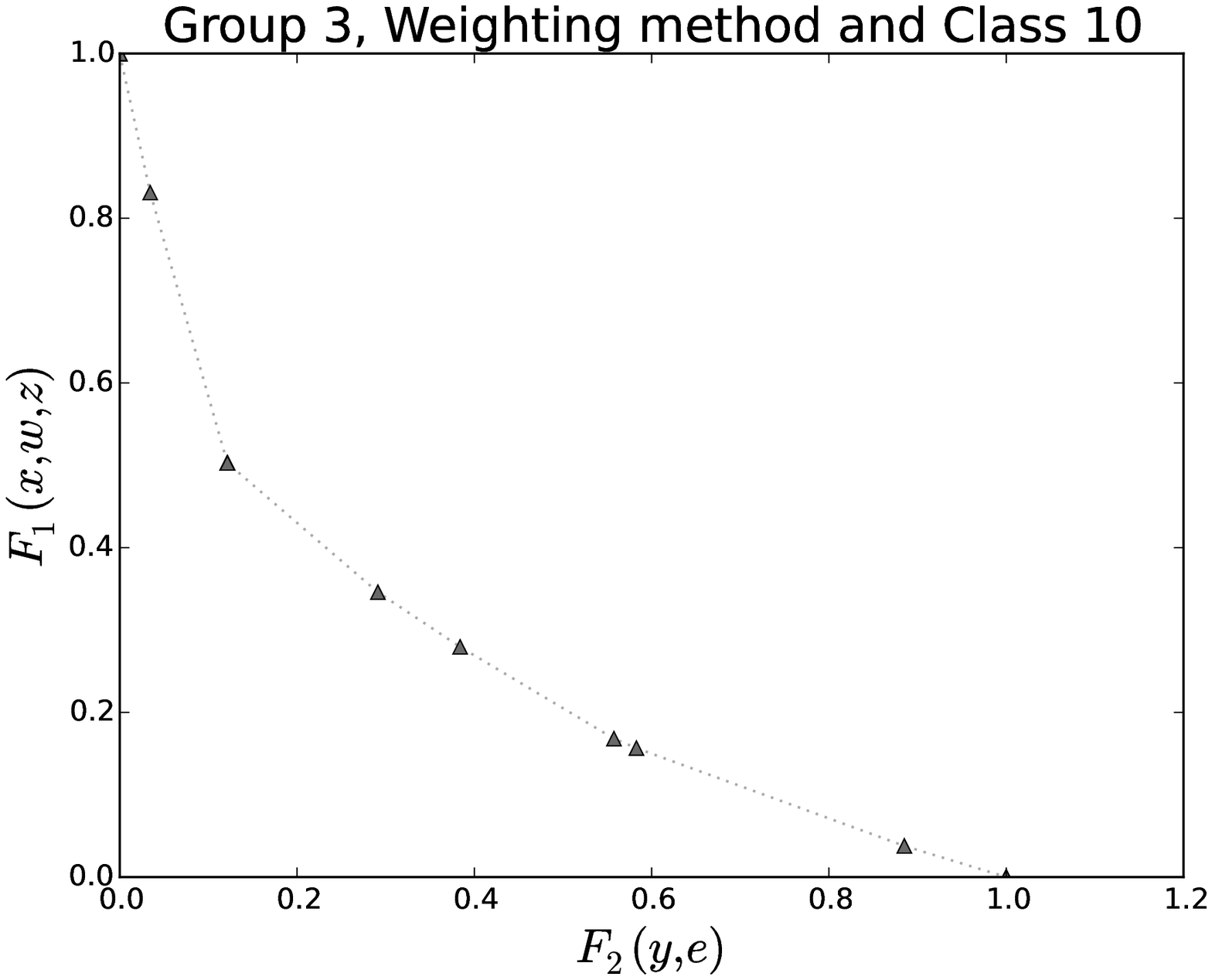}
\caption{Pareto fronts from the weighting and $\varepsilon$-constraint method}\label{fig:Figura2}
\end{center}
\end{figure}

In Figure~\ref{fig:curvaParetoGrupo3Erestrito}, we decided to illustrate the Pareto front obtained from the $\varepsilon$-constraint method for one instance of Group 3 in Class 12, because of the instances with significant dimension in Group 3 that present plenty of points in Pareto front. Each point in the chart represents LSP total costs in opposition to CSP total costs by varying the bounds $\varepsilon_2$. Note that as CSP costs increase LSP costs decrease, and vice versa, evidencing the trade-off between them. Thus we can note that there are many solutions for the bi-objective integrated problem. Furthermore, this behavior remains in all the test groups and instance classes even any other dimension of the instance. Figure~\ref{fig:Figura2} illustrates some results obtained by $\varepsilon$-constraint method and weighting method for Classes 3, 5 and 10.

\begin{figure}[t]
\begin{center}
\includegraphics[width=5cm]{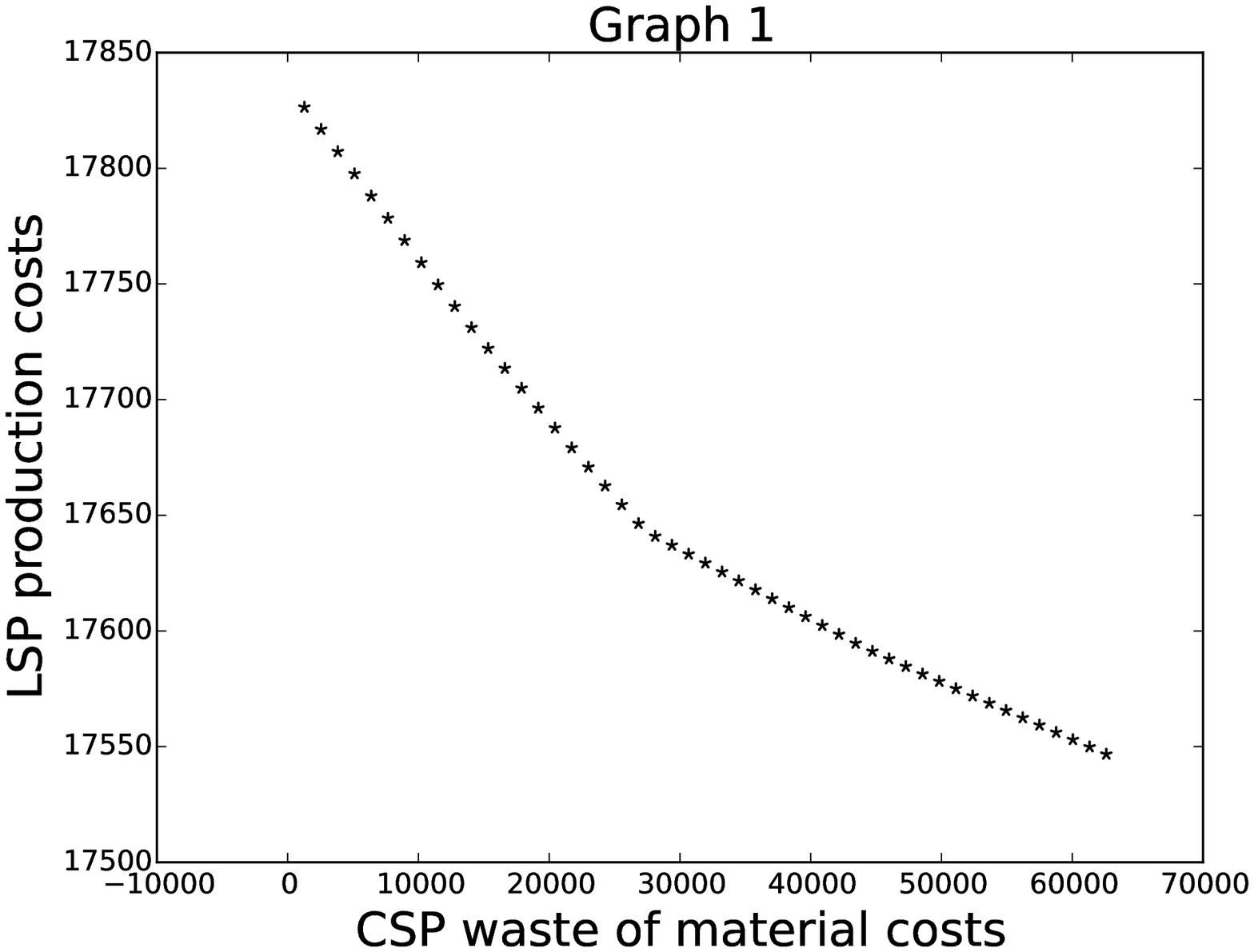} 
\includegraphics[width=5cm]{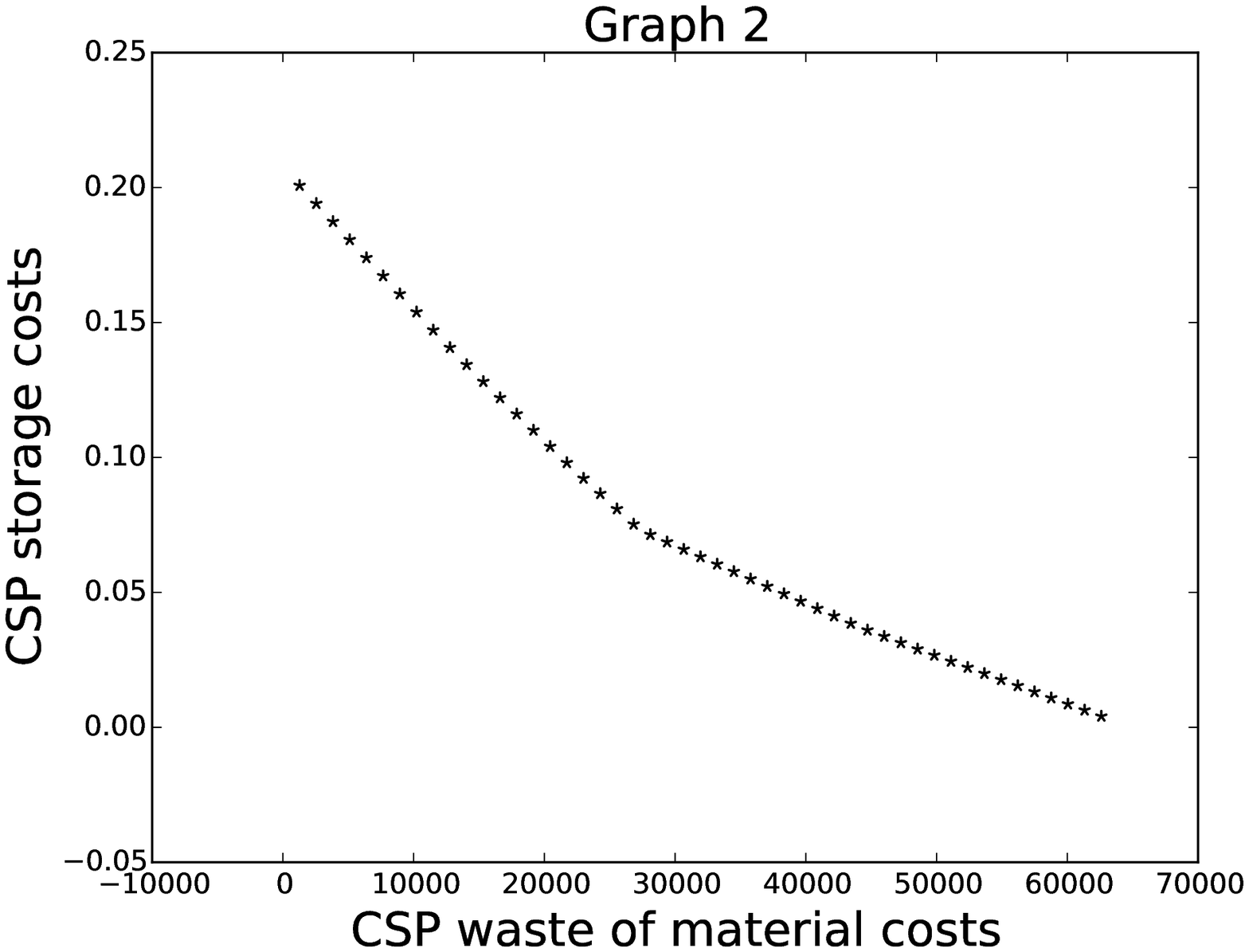} 
\includegraphics[width=5cm]{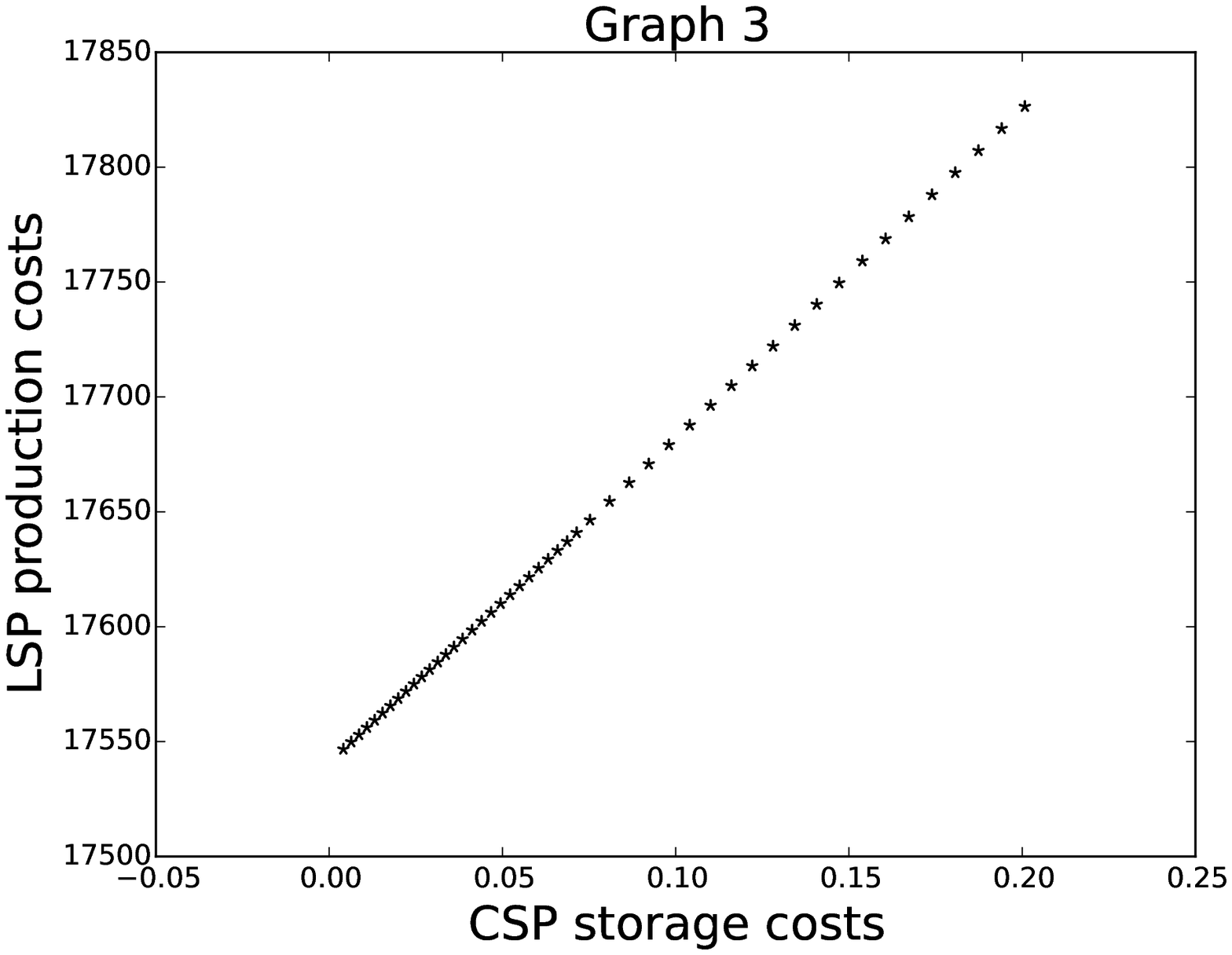}
\includegraphics[width=5cm]{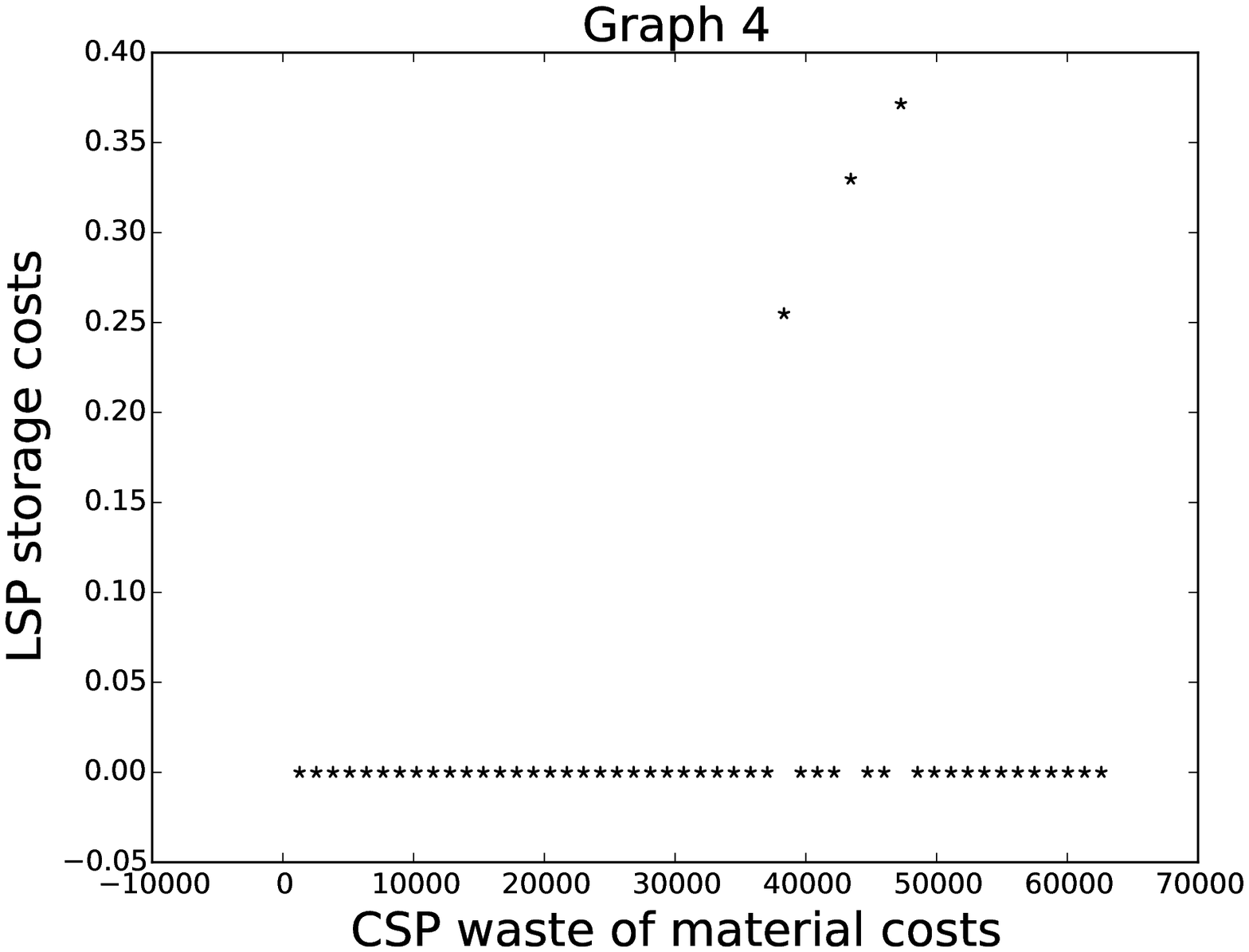}
\includegraphics[width=5cm]{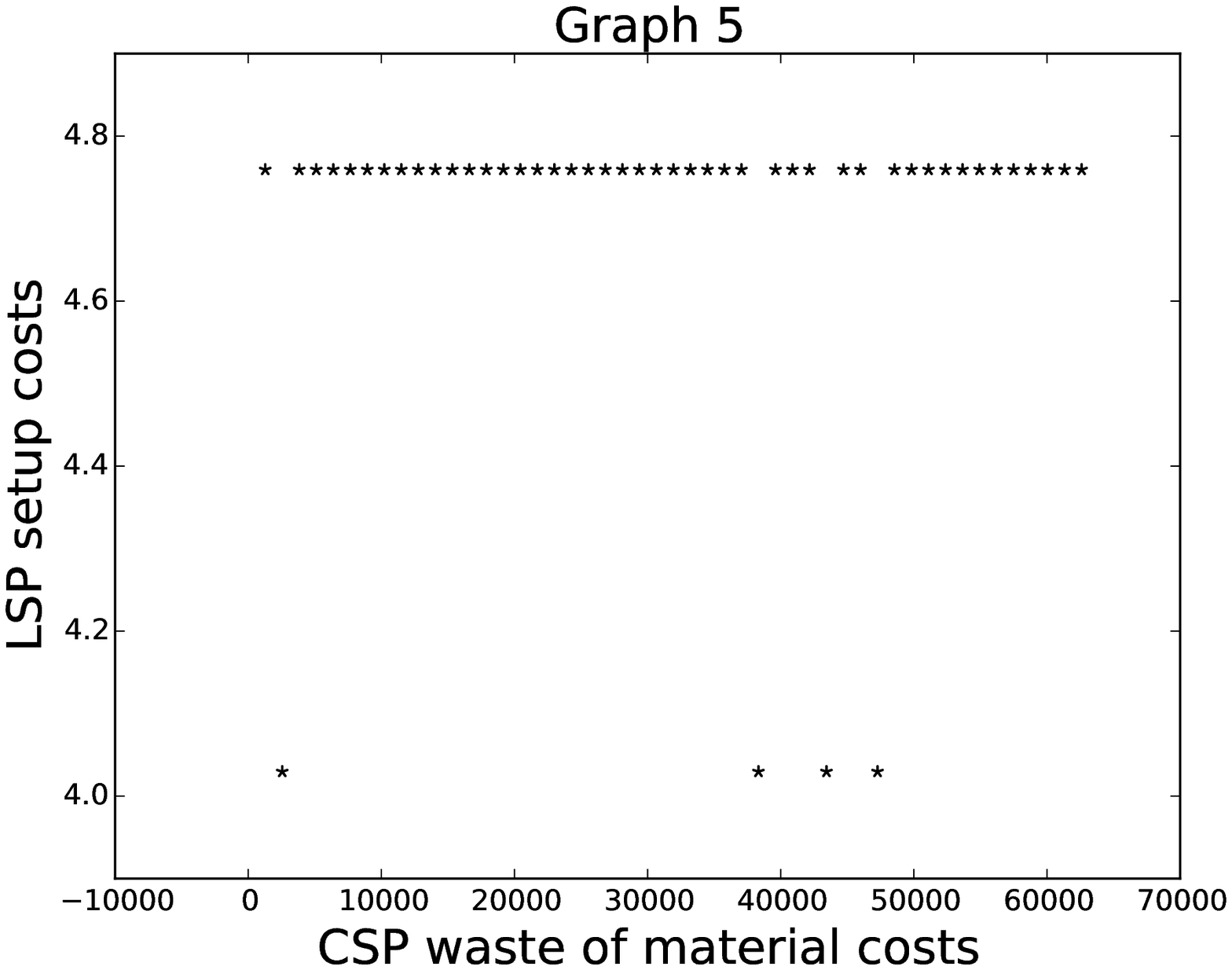}
\includegraphics[width=5cm]{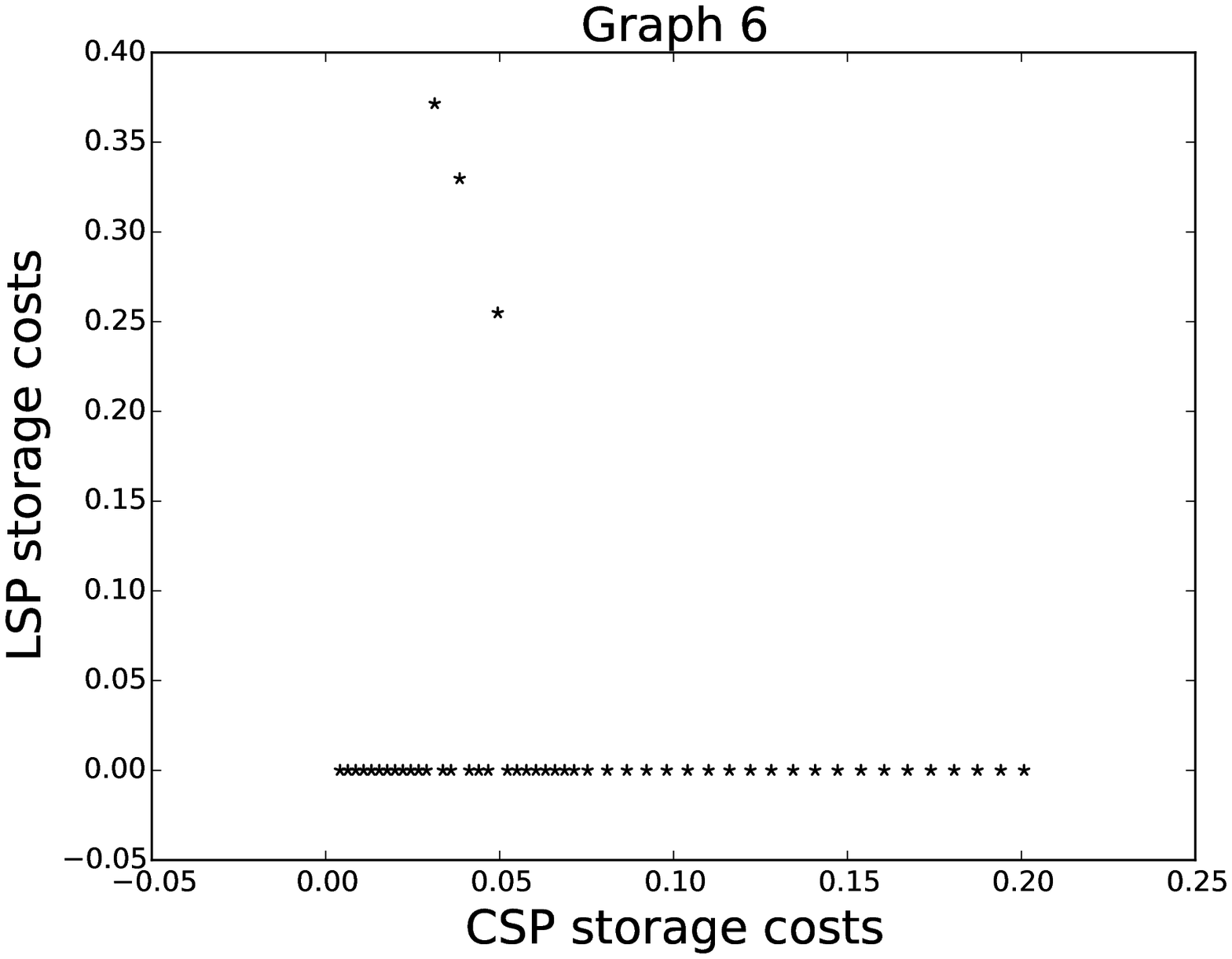}
\includegraphics[width=5cm]{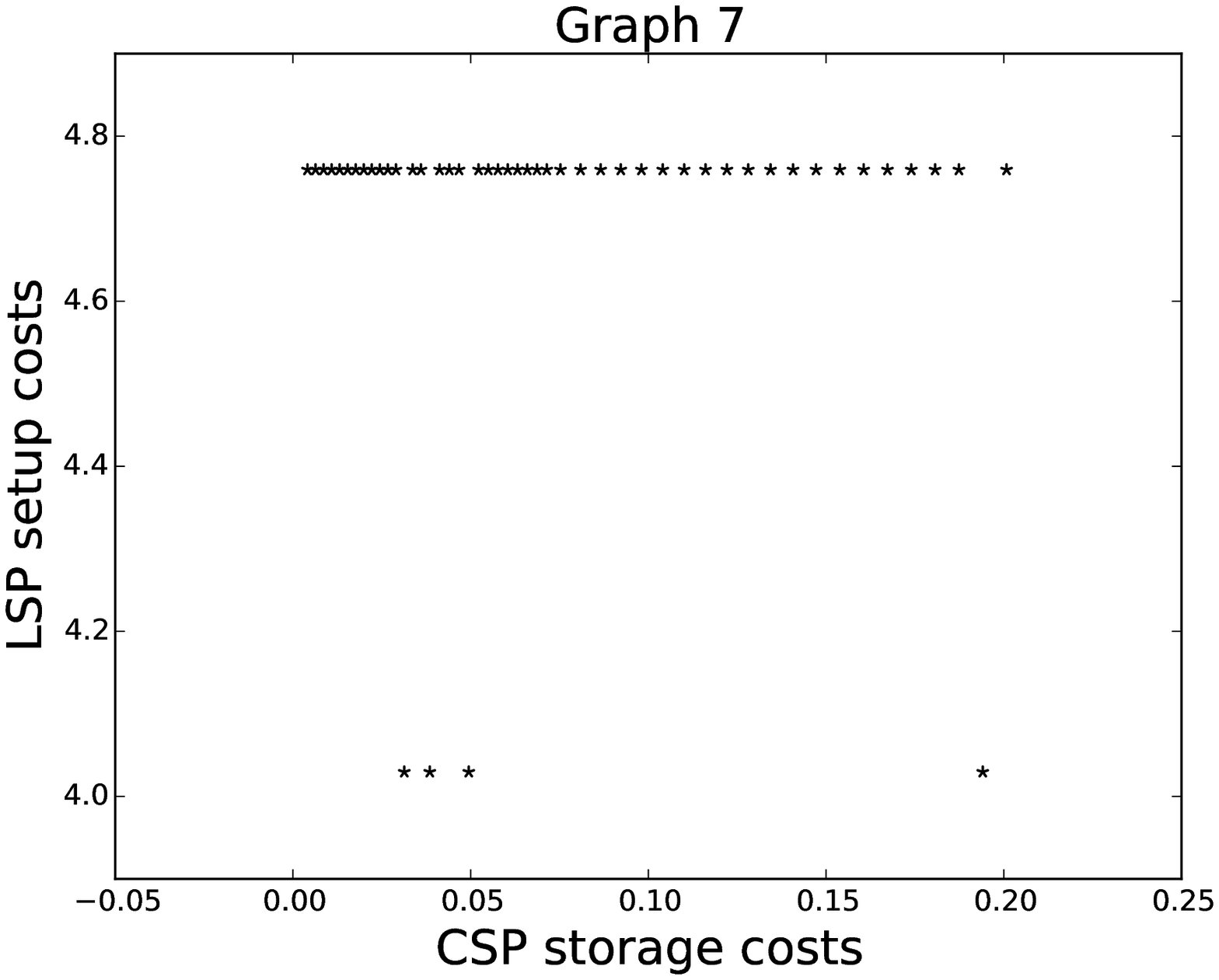}
\includegraphics[width=5cm]{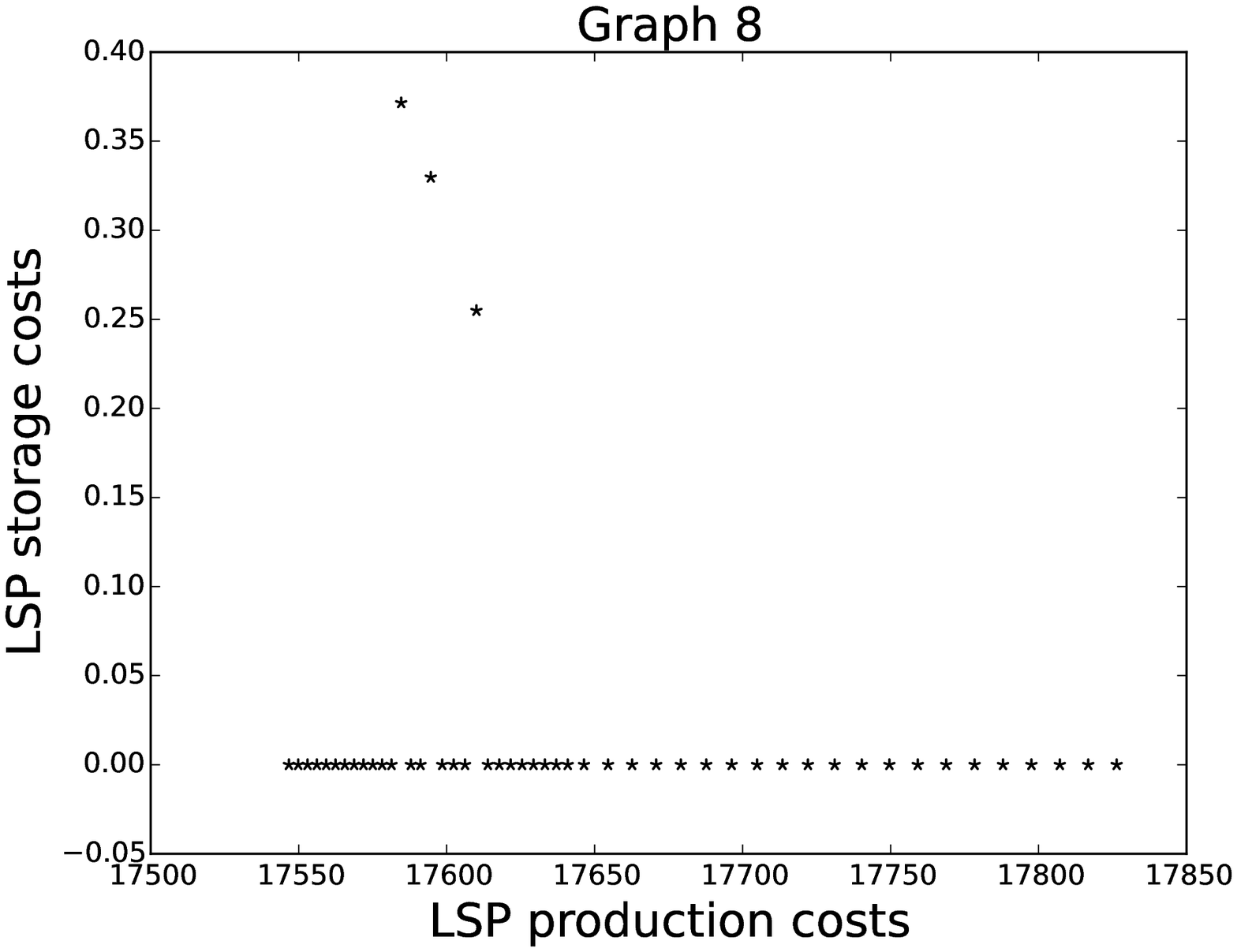}
\includegraphics[width=5cm]{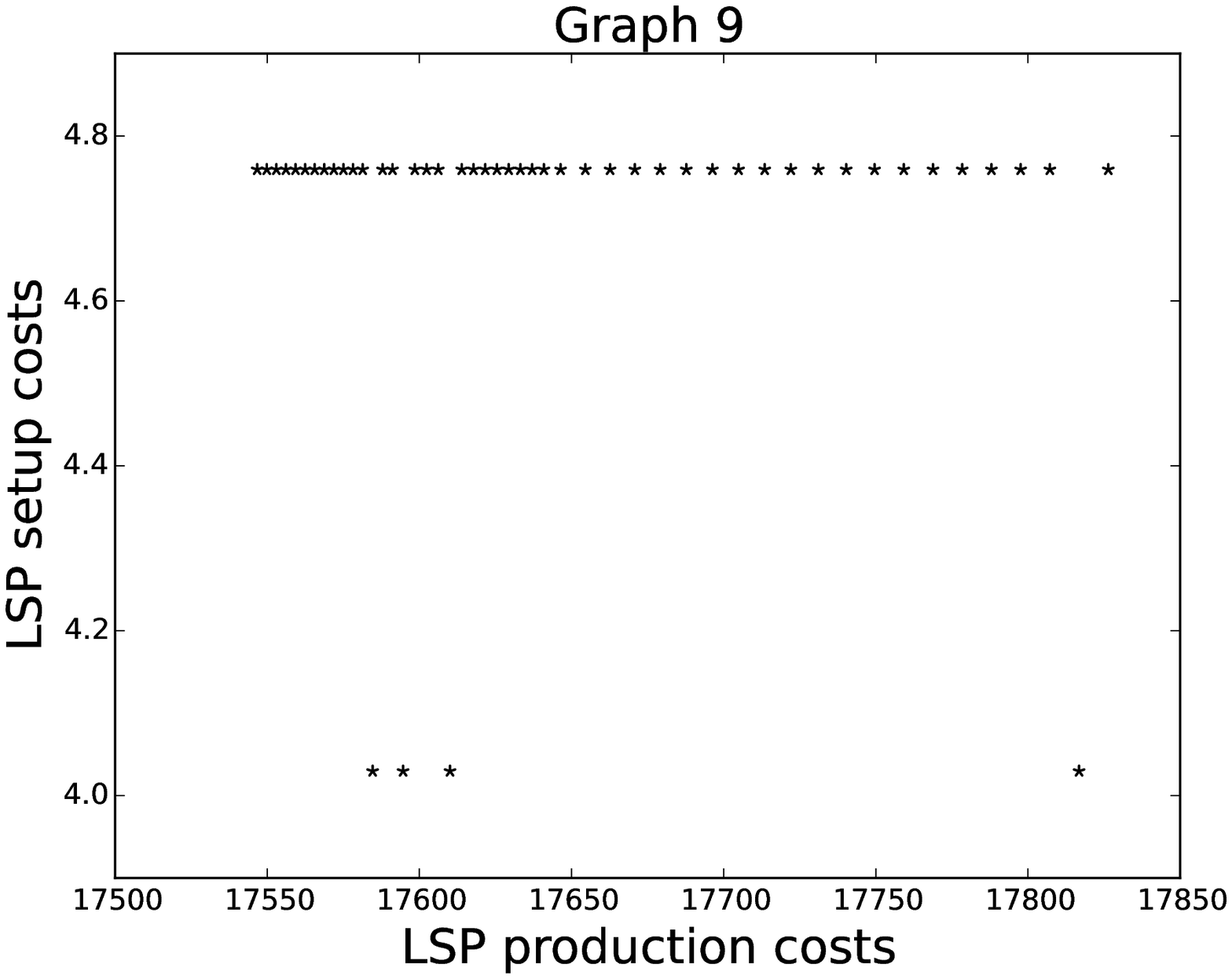}
\includegraphics[width=5cm]{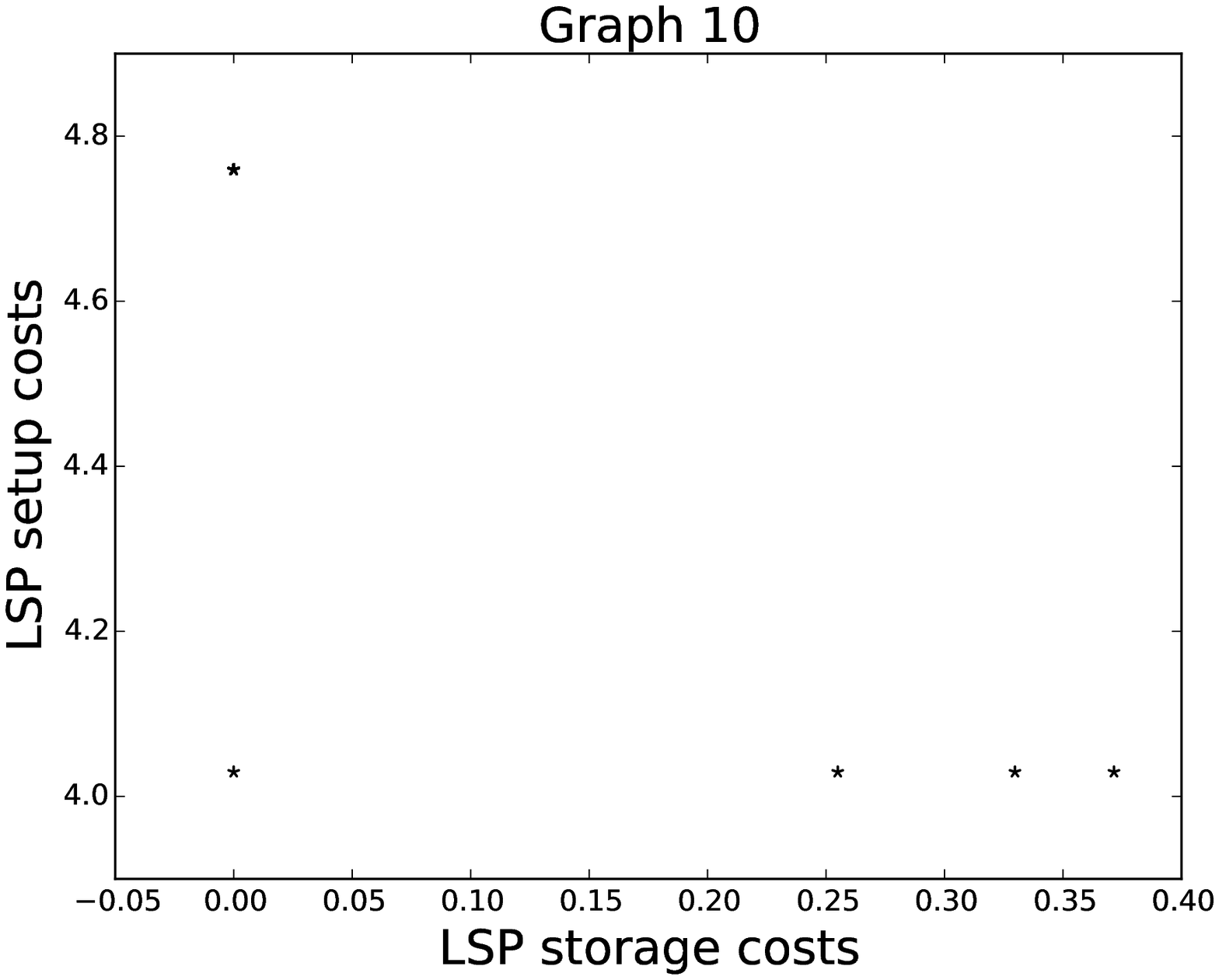}
\caption{Comparisons from the dispersion graphs}\label{fig:dispersao}
\end{center}
\end{figure}

The trade-off noted in the results indicates a negative correlation between the objective functions. We think this is important to make an additional detailed analysis on each variable, reflected in various costs linked to the production process, to detect its interdependence.

\begin{figure}[t]
\begin{center}
\includegraphics[width=5cm]{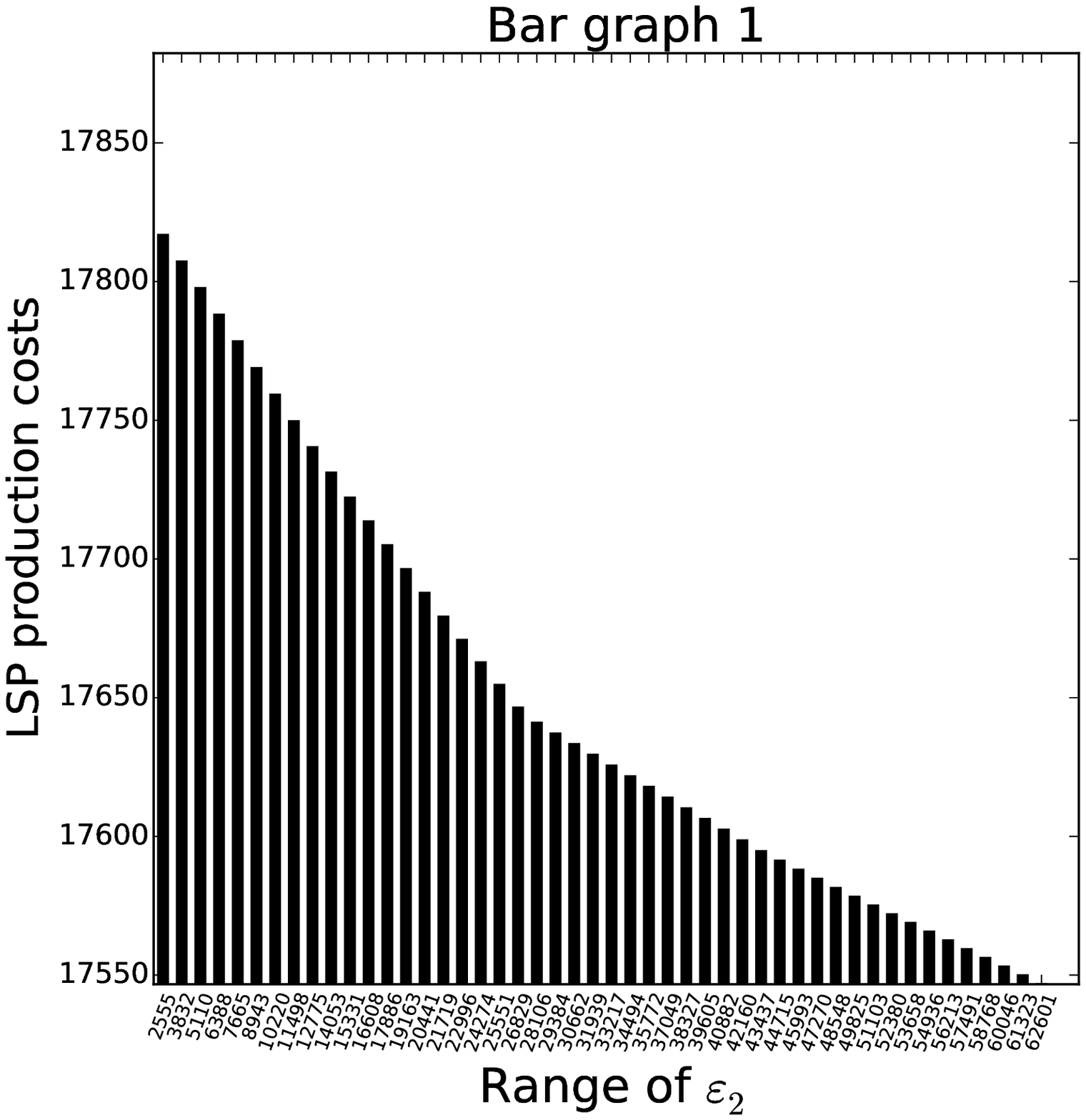} 
\includegraphics[width=5cm]{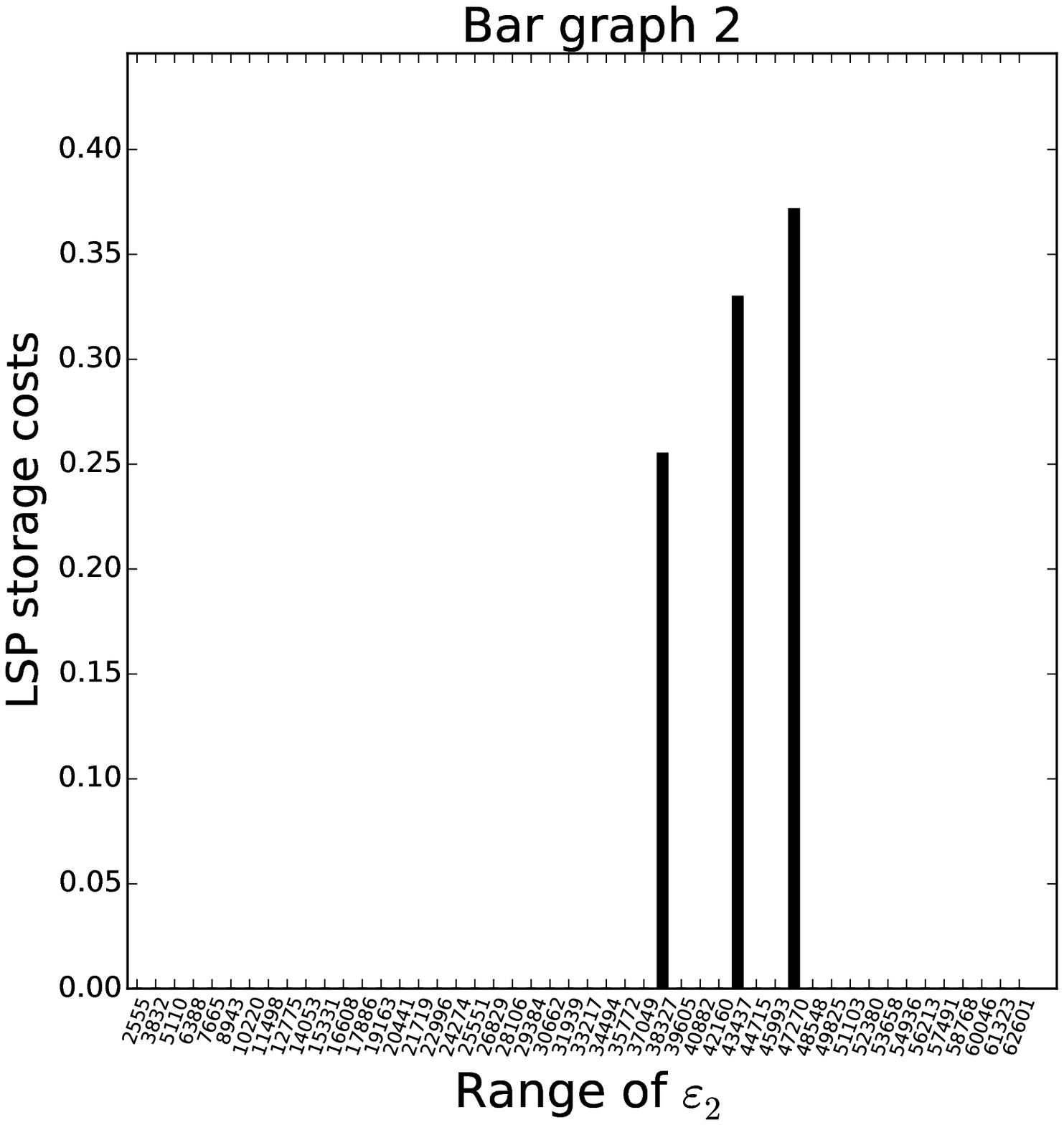} 
\includegraphics[width=5cm]{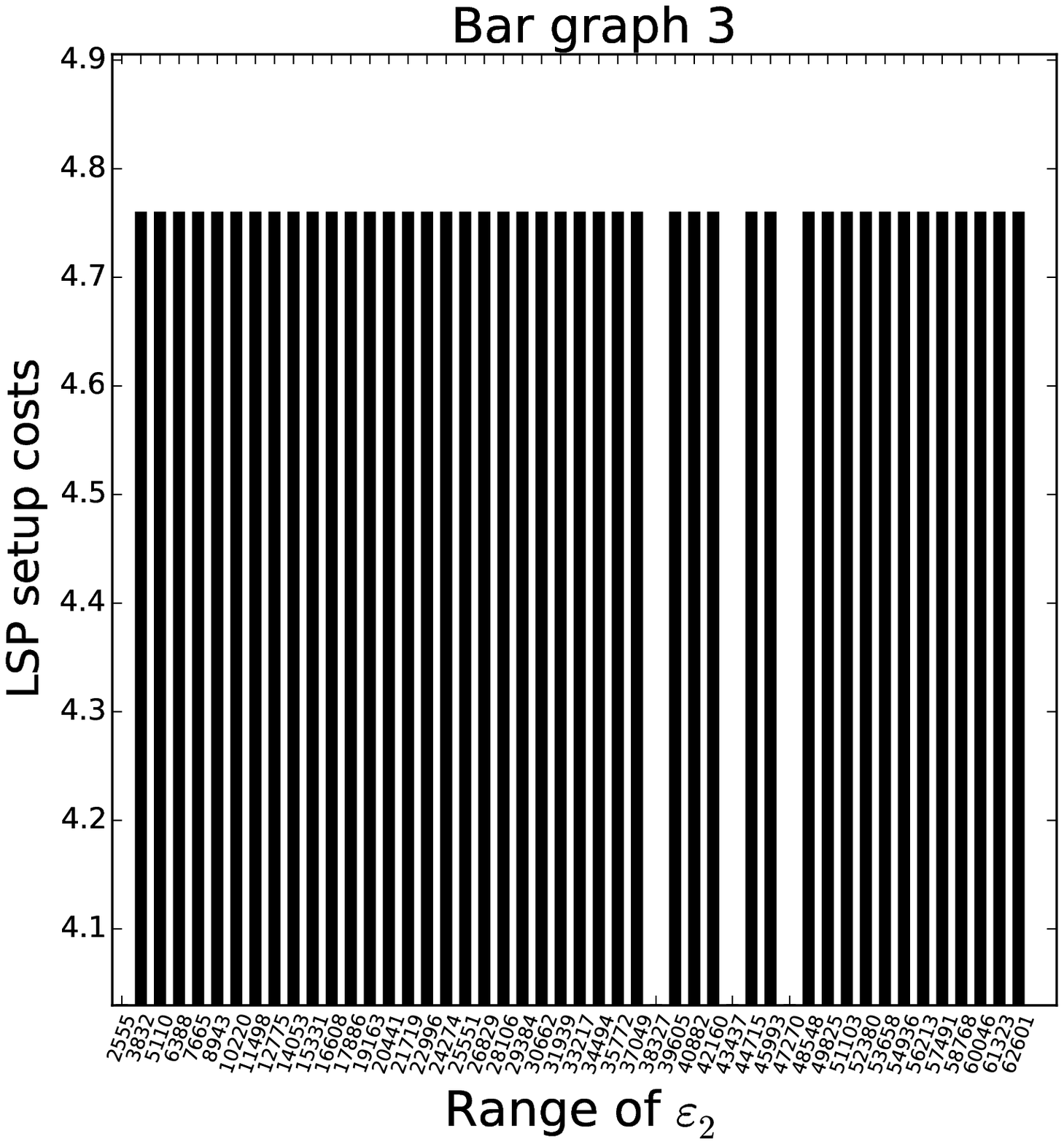}
\includegraphics[width=5cm]{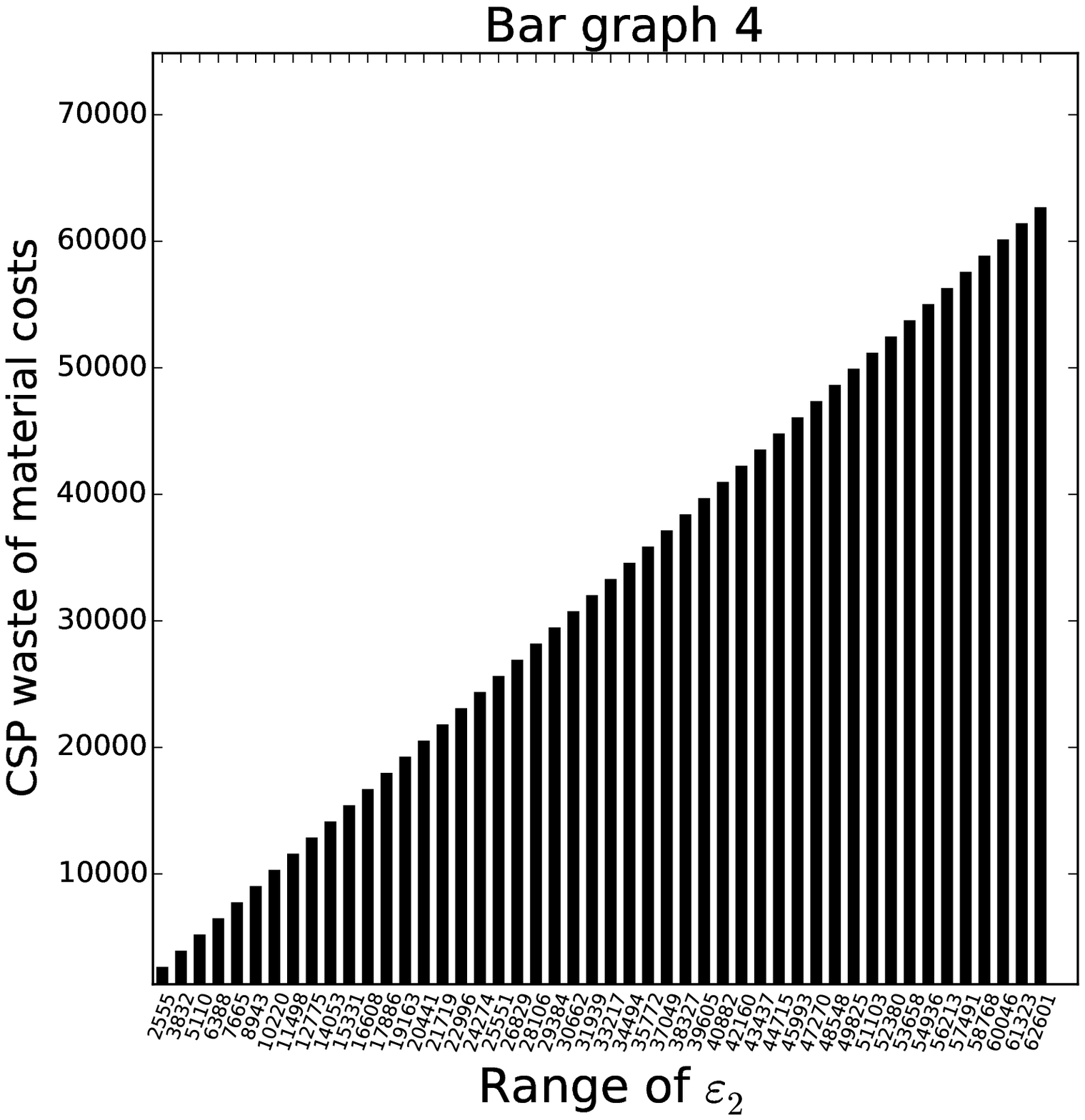}
\includegraphics[width=5cm]{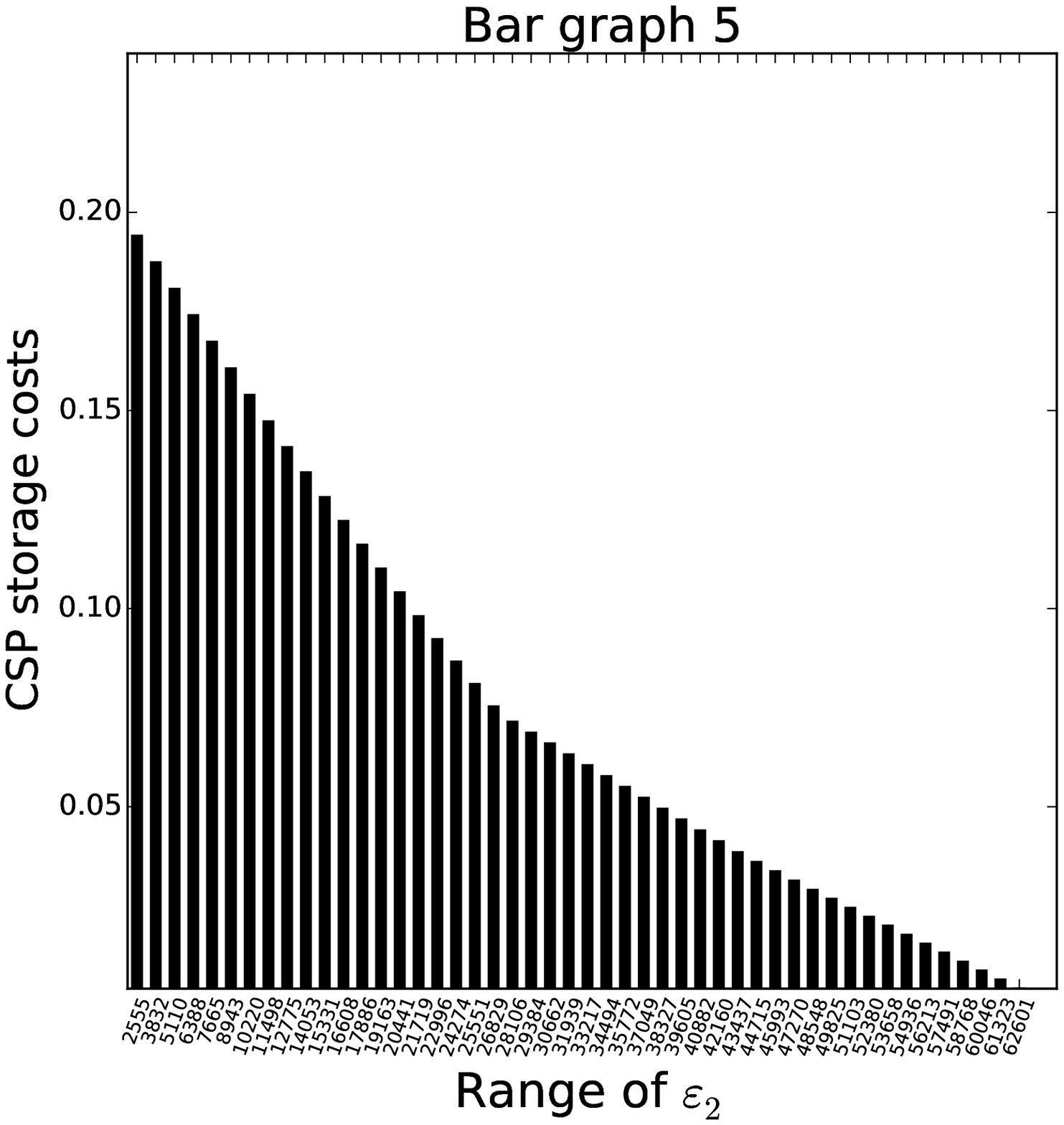}
\includegraphics[width=5cm]{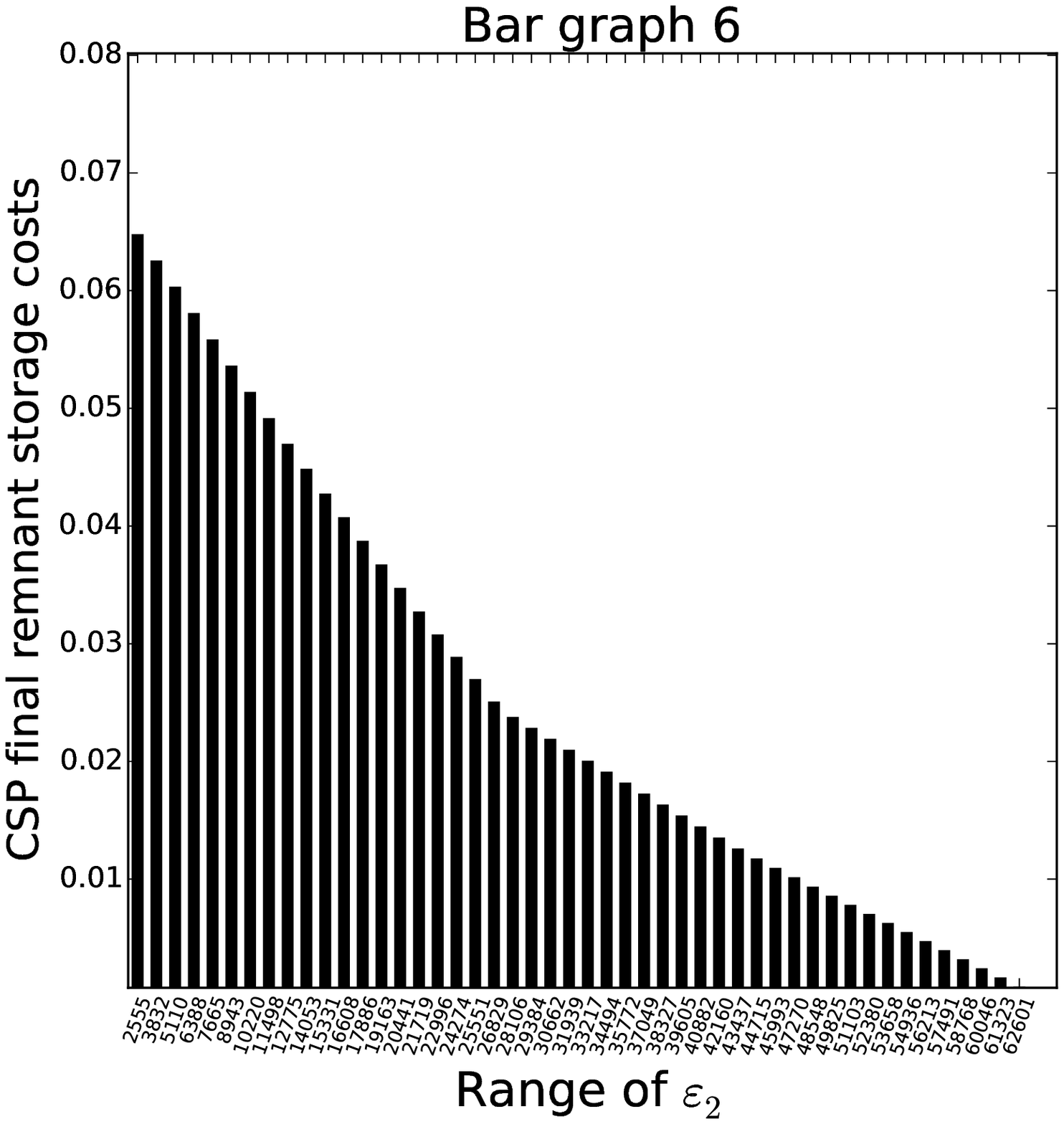}
\caption{Comparisons from the bar graphs and variations of $\varepsilon_2$}\label{fig:graficodebarra}
\end{center}
\end{figure}

Figures~\ref{fig:dispersao} and~\ref{fig:graficodebarra} show respectively the dispersion graphs and bar graphs obtained from $\varepsilon$-constraint method for an instance of Group 3 in Class 11. They illustrate several tests conducted by splitting the bi-objective function in other suitable functions to compare each other. Namely, 

\noindent{$g_1(x)\equiv\sum\limits_{t=1}^{T} \sum\limits_{m=1}^{M} \sum\limits_{k=1}^{K} (c_{kmt} x_{kmt})$} is LSP production costs; 

\noindent{$g_2(w)\equiv\sum\limits_{t=1}^{T} \sum\limits_{m=1}^{M} \sum\limits_{k=1}^{K} (h_{kt} b_{km} w_{kmt})$} is LSP storage costs;

\noindent{$g_3(z)\equiv\sum\limits_{t=1}^{T} \sum\limits_{m=1}^{M} \sum\limits_{k=1}^{K} (s_{kmt} z_{kmt})$} is LSP setup costs; 

\noindent{$g_4(y)\equiv\sum\limits_{t=1}^{T} \sum\limits_{k=1}^{K} cp_{kt} \sum\limits_{m=1}^{M} \sum\limits_{j=1}^{N_m} (p_{jm} y^{j}_{kmt})$} is CSP waste of material costs; and 

\noindent{$g_5(e)\equiv\sum\limits_{t=1}^{T} \sum\limits_{k=1}^{K} \sum\limits_{i \in S(k)} \sigma_{it} \eta_{ik} e_{ikt}$} is CSP storage costs. 

\noindent Each instance is solved using $\varepsilon$-constraint method and we evaluate each one of these split functions at the Pareto optimal points found, and we compared the split functions each other. 

By varying the bound $\varepsilon_2$ in its interval range, Figures~\ref{fig:dispersao} and~\ref{fig:graficodebarra} seek to represent if there is interdependence/correlation between two split functions. Note in Graph 1 of Figure~\ref{fig:dispersao} the negative correlation between LSP production and CSP waste of material costs. In Graph 2, the same happens between CSP storage and CSP waste of material costs. However, in Graph 3, the correlation between LSP production and CSP storage costs is positive.

The other graphs in Figure~\ref{fig:dispersao} give no evidence of correlations between other split functions. The same behavior remains in all test groups and instance classes. Therefore, through the above results, we can conclude that increasing the amount of produced objects implies reducing the waste of material. There is also an inner negative correlation in CSP, i.e., the storage costs increase when the waste of material decreases.

From the bar graphs in Figure~\ref{fig:graficodebarra}, we have a more precise analysis of the variables behavior (measured by costs) by changing $\varepsilon_2$. Unlike in previous split functions, we include CSP final remnant storage costs, measured by {\small$\sum_{k=1}^{K} \sum_{i \in S(k)} \sigma_{iT} \eta_{ik} e_{ikT}$}. That means we produce additional pieces to the required demands and store them at the end of the last period $T$. We can note in Bar graphs 1, 5 and 6 that, as we gradually raise the bound $\varepsilon_2$, the compared split function reduce. On the other hand, in Bar graph 4, we note an increase in CSP waste of material costs. In contrast, in Bar graphs 2 and 3, it is not possible to establish any conflicted relation.

\begin{table}[t]
\caption{Pearson product-moment correlation coefficient among the split functions}
{\footnotesize \begin{tabular}{cccccccccccc}
\toprule
\multicolumn{1}{c}{} & \multicolumn{1}{c}{$g_1(x) \times$} & \multicolumn{1}{c}{$g_4(y) \times$} & \multicolumn{1}{c}{$g_1(x) \times$ } & \multicolumn{1}{c}{$g_2(w) \times$} & \multicolumn{1}{c}{$g_3(z) \times$} & \multicolumn{1}{c}{$g_2(w) \times$} & \multicolumn{1}{c}{$g_3(z) \times$} & \multicolumn{1}{c}{$g_1(x) \times$} & \multicolumn{1}{c}{$g_1(x)\times$} & \multicolumn{1}{c}{$g_2(w) \times$} & \multicolumn{1}{c}{$F_1\times$} \\ 
\multicolumn{1}{c}{} & \multicolumn{1}{c}{$g_4(y)$} & \multicolumn{1}{c}{$g_5(e)$} & \multicolumn{1}{c}{$g_5(e)$} & \multicolumn{1}{c}{$g_4(y)$} & \multicolumn{1}{c}{$g_4(y)$} & \multicolumn{1}{c}{$g_5(e)$} & \multicolumn{1}{c}{$g_5(e)$} & \multicolumn{1}{c}{$g_2(w)$} & \multicolumn{1}{c}{$g_3(z)$} & \multicolumn{1}{c}{$g_3(z)$} & \multicolumn{1}{c}{$F_2$} \\\midrule\midrule
\multicolumn{1}{c}{{\bf Class}} & \multicolumn{11}{c}{{\bf Group 1}} \\ \cmidrule{2-12}
1 &\cellcolor{gray!30} -0.90 &\cellcolor{gray!30} -0.85 & 0.92 & 0.44 & -0.09 & -0.41 & 0.09 & -0.44 & 0.13 & -0.43 & \cellcolor{gray!30} -0.89 \\ 
2 &\cellcolor{gray!30} -0.91 &\cellcolor{gray!30} -0.85 & 0.91 & 0.48 & -0.05 & -0.45 & 0.09 & -0.46 & 0.10 & -0.46 & \cellcolor{gray!30} -0.91 \\ 
3 &\cellcolor{gray!30} -0.82 & \cellcolor{gray!30} -0.83 & 0.99 & 0.33 & 0.08 & -0.34 & -0.12 & -0.34 & -0.12 & -0.38 & \cellcolor{gray!30} -0.82 \\ 
4 &\cellcolor{gray!30} -0.86 &\cellcolor{gray!30} -0.84 & 0.92 & 0.24 & 0.13 & -0.21 & -0.14 & -0.24 & -0.13 & -0.45 & \cellcolor{gray!30} -0.86 \\ 
5 &\cellcolor{gray!30} -0.91 &\cellcolor{gray!30} -0.88 & 0.95 & 0.33 & 0.08 & -0.30 & -0.05 & -0.33 & -0.05 & -0.42 & \cellcolor{gray!30} -0.91 \\ 
6 &\cellcolor{gray!30} -0.84 & -0.75 & 0.81 & 0.26 & 0.15 & -0.25 & -0.09 & -0.29 & -0.14 & -0.35 & -0.71 \\ 
7 &\cellcolor{gray!30} -0.89 &\cellcolor{gray!30} -0.82 & 0.93 & 0.22 & 0.04 & -0.24 & -0.06 & -0.25 & -0.03 & -0.37 & \cellcolor{gray!30} -0.89 \\ 
8 &\cellcolor{gray!30} -0.83 &\cellcolor{gray!30} -0.81 & 0.87 & 0.04 & 0.02 & -0.08 & 0.00 & -0.09 & 0.01 & -0.41 & \cellcolor{gray!30} -0.82 \\ \cmidrule{2-12}
\multicolumn{1}{c}{{\bf Average}} &\cellcolor{gray!30} -0.87 &\cellcolor{gray!30} -0.83 & 0.91 & 0.29 & 0.04 & -0.29 & -0.04 & -0.30 & -0.03 & -0.41 & \cellcolor{gray!30} -0.85 \\ \midrule\midrule
\multicolumn{1}{c}{{\bf Class}} & \multicolumn{11}{c}{{\bf Group 2}} \\ \cmidrule{2-12}
1 &\cellcolor{gray!30} -0.93 &\cellcolor{gray!30} -0.93 & 0.99 & 0.44 & -0.01 & -0.45 & 0.01 & -0.45 & 0.01 & -0.36 &\cellcolor{gray!30} -0.93 \\ 
2&\cellcolor{gray!30} -0.90&\cellcolor{gray!30} -0.92 &0.98&0.39&0.02&-0.40&0.02&-0.39&0.03&-0.55&\cellcolor{gray!30} -0.92\\
3 &\cellcolor{gray!30} -0.85 &\cellcolor{gray!30} -0.82 & 0.94 & 0.34 & 0.06 & -0.32 & -0.09 & -0.35 & -0.05 & -0.42 &\cellcolor{gray!30} -0.85 \\ 
4 &\cellcolor{gray!30} -0.87 & \cellcolor{gray!30} -0.87 & 0.97 & 0.18 & 0.19 & -0.19 & -0.19 & -0.22 & -0.17 & -0.40 & \cellcolor{gray!30} -0.87  \\\cmidrule{2-12}
\multicolumn{1}{c}{{\bf Average}} & \cellcolor{gray!30} -0.89 & \cellcolor{gray!30} -0.89 & 0.97 & 0.39 & 0.02 & -0.39 & -0.02 & -0.40 & 0.00 & -0.44 &\cellcolor{gray!30} -0.90 \\ \midrule\midrule
\multicolumn{1}{c}{{\bf Class}} & \multicolumn{11}{c}{{\bf Group 3}} \\ \cmidrule{2-12}
1 &\cellcolor{gray!30} -0.89 &\cellcolor{gray!30} -0.83 & 0.93 & 0.24 & -0.01 & -0.23 & 0.05 & -0.26 & 0.11 & -0.41 &\cellcolor{gray!30} -0.91 \\
2 &\cellcolor{gray!30} -0.87 &\cellcolor{gray!30} -0.82 & 0.88 & 0.36 & 0.13 & -0.39 & -0.02 & -0.35 & 0.00 & -0.52 &\cellcolor{gray!30} -0.89 \\ 
3 &\cellcolor{gray!30} -0.93 &\cellcolor{gray!30} -0.92 & 0.99 & 0.23 & 0.14 & -0.21 & -0.10 & -0.23 & -0.07 & -0.37 &\cellcolor{gray!30} -0.93 \\ 
4 &\cellcolor{gray!30} -0.83 &\cellcolor{gray!30} -0.84 & 0.98 & 0.22 & 0.23 & -0.19 & -0.21 & -0.19 & -0.18 & -0.39 &\cellcolor{gray!30} -0.83 \\ 
5 &\cellcolor{gray!30} -0.87 &\cellcolor{gray!30} -0.88 & 1.00 & 0.14 & 0.20 & -0.14 & -0.21 & -0.14 & -0.20 & -0.44 &\cellcolor{gray!30} -0.87 \\ 
6 &\cellcolor{gray!30} -0.90 &\cellcolor{gray!30} -0.90 & 0.99 & 0.35 & 0.25 & -0.26 & -0.28 & -0.27 & -0.25 & -0.43 &\cellcolor{gray!30} -0.90 \\ 
7 &\cellcolor{gray!30} -0.88 &\cellcolor{gray!30} -0.89 & 0.99 & 0.05 & 0.27 & -0.08 & -0.26 & -0.10 & -0.24 & -0.25 &\cellcolor{gray!30} -0.88 \\ 
8 &\cellcolor{gray!30} -0.87 &\cellcolor{gray!30} -0.88 & 0.99 & 0.12 & 0.19 & -0.09 & -0.24 & -0.11 & -0.22 & -0.50 &\cellcolor{gray!30} -0.87 \\ 
9 &\cellcolor{gray!30} -0.92 &\cellcolor{gray!30} -0.93 & 0.99 & 0.00 & 0.19 & 0.00 & -0.26 & -0.01 & -0.25 & -0.53 &\cellcolor{gray!30} -0.92 \\ 
10 &\cellcolor{gray!30} -0.91 &\cellcolor{gray!30} -0.91 & 1.00 & 0.09 & 0.25 & -0.10 & -0.25 & -0.10 & -0.24 & -0.59 &\cellcolor{gray!30} -0.91 \\ 
11 &\cellcolor{gray!30} -0.81 &\cellcolor{gray!30} -0.85 & 0.96 & -0.11 & 0.07 & 0.10 & -0.11 & 0.10 & -0.14 & -0.47 &\cellcolor{gray!30} -0.89 \\ 
12 &\cellcolor{gray!30} -0.83 &\cellcolor{gray!30} -0.83 & 0.97 & 0.07 & 0.18 & -0.01 & -0.22 & -0.02 & -0.21 & -0.69 &\cellcolor{gray!30} -0.84 \\\cmidrule{2-12} 
\multicolumn{1}{c}{{\bf Average}} &\cellcolor{gray!30} -0.87 &\cellcolor{gray!30} -0.87 & 0.97 & 0.15 & 0.18 & -0.13 & -0.18 & -0.14 & -0.16 & -0.47 &\cellcolor{gray!30} -0.89 \\ \bottomrule
\end{tabular}}
\label{tab:correlation}
\end{table}

\begin{table}[t]
\caption{Number of different Pareto solutions, running time, number of variables and number of constraints}
\footnotesize{\begin{tabular}{lcrrrrrrrrrrrr}
\toprule
& & \multicolumn{4}{c}{\bf Group 1} & \multicolumn{4}{c}{\bf Group 2} & \multicolumn{4}{c}{\bf Group 3} \\\midrule
\bf Method & \bf Class & \it nd & time (s) & \it nv & \it nc &\it nd & time (s) & \it nv &\it nc & \it nd & time (s) & \it nv & \it nc \\ \midrule
& 1 & 39.9 & 24.6 & 110.6 & 31.0 & 41.9 & 136.2 & 120.2 & 31.0 & 34.2 & 5.4 & 121.5 & 31.0 \\ 
& 2 & 43.6 & 141.0 & 147.0 & 41.0 & 44.0 & 150.0 & 152.8 & 41.0 & 38.7 & 3.6 & 185.6 & 41.0 \\ 
& 3 & 42.1 & 117.0 & 120.0 & 34.0 & 47.0 & 36.0 & 259.5 & 34.0 & 48.0 & 4.8 & 307.2 & 34.0 \\ 
& 4 & 44.3 & 170.4 & 160.0 & 45.0 & 38.9 & 149.9 & 309.0 & 45.0 & 47.7 & 7.8 & 407.4 & 45.0 \\ 
& 5 & 43.8 & 136.2 & 123.0 & 37.0 &  &  &  &  & 45.6 & 6.0 & 539.4 & 37.0 \\ 
{\bf $\varepsilon$-Constraint} & 6 & 39.5 & 144.6 & 164.0 & 49.0 &  &  &  &  & 48.0 & 9.6 & 881.6 & 49.0 \\ 
& 7 & 41.7 & 79.2 & 126.0 & 40.0 &  &  &  &  & 45.6 & 10.2 & 1223.7 & 40.0 \\ 
& 8 & 46.0 & 90.0 & 168.0 & 53.0 &  &  &  &  & 48.0 & 14.4 & 1559.6 & 53.0 \\ 
& 9 & & &  & &  &  &  &  & 46.6 & 16.2 & 1892.4 & 43.0 \\ 
& 10 & & &  & &  &  &  &  & 45.6 & 25.2 & 3414.0 & 57.0 \\ 
& 11 & & &  & &  &  &  &  & 43.4 & 27.0 & 3430.8 & 46.0 \\ 
& 12 & & &  & &  &  &  &  & 45.2 & 39.0 & 4909.8 & 61.0 \\ \midrule\midrule
& 1 & 5.5 & 11.4 & 110.6 & 30.0 & 5.7 & 25.8 & 120.2 & 30.0 & 2.2 & 2.4 & 121.5 & 30.0 \\ 
& 2 & 4.7 & 47.4 & 147.0 & 40.0 & 4.4 & 43.8 & 152.8 & 40.0 & 3.2 & 1.8 & 185.6 & 40.0 \\ 
& 3 & 6.1 & 62.4 & 120.0 & 33.0 & 6.0 & 44.4 & 259.5 & 33.0 & 4.1 & 1.8 & 307.2 & 33.0 \\ 
& 4 & 6.8 & 76.2 & 160.0 & 44.0 & 4.9 & 48.7 & 309.0 & 44.0 & 6.2 & 3.0 & 407.4 & 44.0 \\ 
& 5 & 7.0 & 54.6 & 123.0 & 36.0 &  &  &  &  & 5.4 & 2.4 & 539.4 & 36.0 \\ 
{\bf Weighting} & 6 & 6.2 & 19.2 & 164.0 & 48.0 &  &  &  &  & 4.8 & 3.0 & 881.6 & 48.0 \\ 
& 7 & 6.2 & 42.6 & 126.0 & 39.0 &  &  & &  & 7.5 & 4.2 & 1223.7 & 39.0 \\ 
& 8 & 6.9 & 56.4 & 168.0 & 52.0 &  &  & &  & 8.6 & 5.4 & 1559.6 & 52.0 \\ 
& 9 & & &  & &  &  & &  & 6.7 & 6.0 & 1892.4 & 42.0 \\ 
& 10 & & & & &  &  & &  & 9.2 & 9.0 & 3414.0 & 56.0 \\ 
& 11 & & & & &  &  & &  & 6.7 & 10.2 & 3430.8 & 45.0 \\ 
& 12 & & & & &  &  & &  & 11.2 & 14.4 & 4909.8 & 60.0 \\ \bottomrule
\end{tabular}}
\label{tab:NumberSolutions}
\end{table}

Table~\ref{tab:correlation} and~\ref{tab:NumberSolutions} summarize the relevant data collected in the numerical experiments with the three groups of tests described in this section. Further details on the correlations among the split functions are presented in Table~\ref{tab:correlation}. For each instance and each one of two methods of resolution, we used the data set obtained by evaluating the split functions at the Pareto solution to calculate the Pearson product-moment correlation coefficient (Pearson's $r$) (\cite{lee1988thirteen}). The entries of Table~\ref{tab:correlation} (correlation matrix) are the average Pearson's $r$ of the 20 instances per group, including the two methods of resolution. At the bottom of each group, we also present the average Pearson's $r$ of the 12 classes of instances. The observed values in Table~\ref{tab:correlation} confirm the conclusions presented above. The table cell containing these entries with Pearson's $r$ from -1 to -0.8 is shaded for emphasizing a negative and significant correlation. For example, in agreement with Pearson's correlation, we note in Group 2 of Class 1 that the correlation between $g_1(x)$ and $g_4(y)$ (i.e., between LSP production and CSP waste of material costs) is negative and significant. It is also on average equal to -0.93. While the correlation between $g_2(w)$ and $g_5(e)$ (i.e., between LSP storage and CSP storage costs) is not significant and is on average equal -0.45.  

In order to investigate whether the Pareto optimal solutions are well spread and spaced from each other, we use the following criteria. Let $u_1=(x_1,y_1,z_1,w_1,e_1)$ and $u_2=(x_2,y_2,z_2,w_2,e_2)$ be two Pareto optimal solutions. If $|F_1(u_1)-F_1(u_2)|\geq 10^{-4}$ and $|F_2(u_1)-F_2(u_2)|\geq 10^{-4}$, we consider that $u_1$ and $u_2$ are different. Table~\ref{tab:NumberSolutions} presents further details on the number of different Pareto solutions found by each method of resolution, the running time (in seconds), and the number of variables and constraints of the instances. For each instance, the methods of resolution are performed fifty times in accordance with fifty upper bounds $\varepsilon_2$ or fifty weighting coefficients $p_i \in(0,1)$. The computed time in Table~\ref{tab:NumberSolutions} is the average time for 20 instances per class. For each instance, we compute the number of variables and constraints. The values {\it nv} and {\it nc} are the average number of variables and constraints, respectively, for 20 instances per class. These data provide some idea of the size and difficulty to solve each instance in each class. From the purposes of this work, the running time obtained by our resolution approaches in each group can be considered satisfactory. Moreover, for each instance, each method of resolution can find up to fifty different Pareto solutions. The value {\it nd} in Table~\ref{tab:NumberSolutions} is the average number of different Pareto solutions found for 20 instances per class. We observe that the $\varepsilon$-constraint method found about 68\% to 96\% of possible solutions, while the weighting method found about 4\% to 22\% of possible solutions. In this requirement, the $\varepsilon$-constraint method is more flexible than weighting method, and it can provide a broader range of solutions, which could be used by the decision maker in a more accurate way. 

We assume, from overall outcomes concerning the characteristics of the LSP integrated with CSP, that our multiobjective approach allowed us to understand the individual problems better. While it is needed to obtain a little bit waste of material, it is mandatory to increase the production of objects. This conflict occurs because, when it is available a high amount of objects, in CSP it is possible to select better cutting patterns that contain little waste of material, even if these selected cutting patterns exceed the demand over a particular period and pay storage costs by including more pieces than necessary.

In opposition, when it is needed to produce few objects, it is mandatory to increase the waste of material. This point can be explained from Constraints~(\ref{eq:corte_itens}) and~(\ref{eq:acoplamento}) in the mathematical model. Concerning stock, in~(\ref{eq:acoplamento}) the total number of produced objects plus the stock balance must be equal to the total number of cut objects, and in~(\ref{eq:corte_itens}) the total number of cut objects following the cutting patterns must meet the demand of pieces. Thus, if few objects are produced, i.e., generating low LSP costs, in CSP the cutting patterns are selected to supply the demand, even if the waste of material increases. This fact means that the possibilities of cutting patterns obtained from the combination of other objects will not be chosen. For clarity, we present one small numerical example. 

Suppose that the demand is 2 pieces of 92cm and 2 pieces of 115cm of length. If only one object of 460cm is available, the cutting pattern to meet this demand must be 2$\times$92 + 2$\times$115 = 414, which generates the waste of material equal to 46cm. Now, suppose that 2 objects of 460cm are available, the cutting patterns that eliminate the waste of material are 5$\times$92 = 460 and 4$\times$115 = 460. Therefore, in the second case, the total waste of material was zero, and 3 pieces of 92cm and 2 pieces of 115cm must be stored. This surplus can be used over the planning horizon.

Each Pareto optimal solution allows the decision maker to decide the production planning better. Therefore, we believe that the multiobjective approach for the integrated problems is relevant, mainly for industries that use multiple criteria ABC analysis (\cite{flores1987implementing}) in its production planning. From a theoretical point of view, our approach has identified that LSP integrated with CSP suggests that it is worthwhile to investigate its possibilities further. For example, the five split functions can be used to create new approaches by considering novel multiobjective models. The correlations identified in this work may presume new 2-, 3- or up to 4-objective optimization models.

\section{Conclusions}\label{sec:conclusions}
We introduce and study a multiobjective mathematical model to deal with the lot sizing problem integrated with cutting stock problem in a paper industry, in which we adequately transformed an earlier mono-objective model into a new bi-objective optimization problem. 

The central aspect of this contribution is to evaluate the trade-off between LSP and CSP by analyzing the cost variations of these two problems when we simultaneously minimize them. Besides, from a suitable manner, we separated the objective function of the multiobjective problem in some other functions to compare, detect, analyze and study potential trade-offs/correlations at the heart of the problem. We illustrate our results through some Pareto front, which highlighted several trade-offs/correlations among the separated functions in our approach. It was clear that increasing LSP total costs imply decreasing CSP total costs, and vice versa. The data confirm many possibilities of the optional Pareto optimal solutions for this integrated problem. 

All correlations among the suitably separated functions are validated and illustrated. The most important ones are the correlations among LSP production, CSP waste of material, and CSP storage costs. Furthermore, we conclude that it is necessary to cut more than the demand to minimize the waste of material, even if this remnant piece is stored in the final period. In fact, we believe that, in practice, the manufacturing industries that have a long planning horizon could use this strategy in its production process.

The conclusions about the trade-off between the two problems can improve the technical decisions in the industrial sector. Mainly because, in opposition to the mono-objective approach, the several solutions obtained from multiobjective approaches provide more opportunities to the decision maker plans the production. We believe that multiobjective approaches are significant for integrating problems, and this could include other ones derived from the industry of furniture, aluminum, etc.

In future work, we plan to study a new multiobjective approach for this integrated problem considering a tri-objective optimization model, by also putting CSP waste of material costs in opposition to CSP storage costs.

\section*{Acknowledgments}
Development agencies partially support the authors. The first author thanks, CAPES (master's program scholarships). The second author is grateful to FAEPEX-Unicamp (grants 519.292-0285/15 and 519.292-885/15). The third author thanks, FAEPEX-Unicamp (master's program scholarships). The fourth author wishes to thank FAPESP (grant 2015/02184-7) and FAEPEX-Unicamp (grant 519.292-262/15).

\nocite{*}

\bibliographystyle{itor}
\bibliography{References}

%\bibliographystyle{itor}
%\bibliography{References}

\end{document}